\documentclass[11pt]{amsart}
\usepackage{amssymb,color}
\usepackage[square,compress,comma, numbers,sort]{natbib}
\usepackage[shortlabels]{enumitem}

\usepackage[colorlinks=true, citecolor=blue, linkcolor=blue]{hyperref}
\usepackage{amsthm}
\usepackage{amsfonts}
\usepackage{mathrsfs}
\usepackage{cleveref}

\definecolor{c30}{rgb}{0.,0.,1.}
\def\Eh#1{{\textcolor{black}{#1}}}
\def\eHe#1{{\textcolor{blue}{#1}}}
\def\eHe#1{{#1}}
\definecolor{c40}{rgb}{0.2,0.2,0.5}
\def\eeh#1{\textcolor{black}{#1}}
\def\tnu#1{\eeh{\overline{T}_{#1}(u)}}

\newcommand{\E}[1]{\mathbb{E}\left\{ #1\right\}}

\newcommand{\pk}[1]{\mathbb{P} \left\{ #1 \right\} }

\def\nvu{\frac 1 {v(u)}}
\def\ndu{\frac 1 {\Delta(u)}}

\def\suk0{\sigma_{\xi_{u,j}}(\vk{0})}

\def\ruk{\rho_{u,j}}
\def\vkt0{(\vk{t},\vk{0})}
\def\vks0{(\vk{s},\vk{0})}

\def\KK{E}

\definecolor{c20}{rgb}{0.,0.7,0.}
\definecolor{c30}{rgb}{0.,0.,1.}
\definecolor{c40}{rgb}{1,0.1,0.7}
\definecolor{c50}{rgb}{1,0,0}
\definecolor{c60}{rgb}{0,0.9,0.1}

\newcommand{\vk}[1]{{#1}}

\newcommand{\ve}{\varepsilon}

\newcommand{\abs}[1]{\left| #1 \right|}
\newcommand{\ABs}[1]{ \biggl \lvert #1 \biggr \rvert}

\newcommand{\R}{\mathbb{R}}

\newcommand{\inr}{\in \R}

\newcommand{\limit}[1]{\lim_{#1 \to   \infty}}

\newcommand{\BQN}{\begin{eqnarray}}
\newcommand{\EQN}{\end{eqnarray}}
\newcommand{\BQNY}{\begin{eqnarray*}}
\newcommand{\EQNY}{\end{eqnarray*}}

\def\eL#1{\textcolor{c50}{#1}}
\def\eL#1{{#1}}

\def\K1#1{\textcolor{cyan}{#1}}
\def\K1#1{#1}
\def\eaE#1{\textcolor{black}{#1}}

\newcommand{\BS}{\begin{sat}}
\newcommand{\ES}{\end{sat}}
\newcommand{\BT}{\begin{theo}}
\newcommand{\ET}{\end{theo}}
\newcommand{\BK}{\begin{korr}}
\newcommand{\EK}{\end{korr}}

\newcommand{\BD}{\begin{de}}
\newcommand{\ED}{\end{de}}

\newcommand{\BIT}{\begin{itemize}}
\newcommand{\EIT}{\end{itemize}}
\newcommand{\BDI}{\begin{description}}
\newcommand{\EDI}{\end{description}}

\newcommand{\BRM}{\begin{remark}}
\newcommand{\ERM}{\end{remark}}

\newcommand{\BEL}{\begin{lem}}
\newcommand{\EEL}{\end{lem}}

\newtheorem{theo}{Theorem}[section]
\newtheorem{sat}[theo]{Proposition}
\newtheorem{de}[theo]{Definition}
\newtheorem{lem}[theo]{Lemma}

\newtheorem{korr}[theo]{Corollary}
\newtheorem{remark}[theo]{Remark}

\newcommand{\prooftheo}[1]{ \textbf{Proof of Theorem} \ref{#1} }
\newcommand{\proofprop}[1]{\textbf{Proof of Proposition} \ref{#1}}
\newcommand{\prooflem}[1]{\textbf{Proof of Lemma} \ref{#1}}

\newcommand{\COM}[1]{}

\newcommand{\QED}{\hfill $\Box$ \\}

\def\vp{\varepsilon}

\def\rw{\rightarrow}

\def\IF{\infty}

\def\LT{\left}
\def\RT{\right}

\def\ehe#1{\textcolor{c50}{#1}}
\def\ehe#1{#1}

\def\dd#1{\textcolor{green}{#1}}
\def\dd#1{#1}
\def\kdd#1{\textcolor{cyan}{#1}}
\def\kdd#1{#1}

\def\I{\mathbb{I}}
\def\rd#1{\textcolor{black}{#1}}

\def\OU{\eHe{\omega(u)}}

\begin{document}

\title
{Sojourns of fractional Brownian motion queues: transient asymptotics}

\author{Krzysztof D\c{e}bicki}
\address{Krzysztof D\c{e}bicki, Mathematical Institute, University of Wroc\l aw, pl. Grunwaldzki 2/4, 50-384 Wroc\l aw, Poland}
\email{Krzysztof.Debicki@math.uni.wroc.pl}

\author{Enkelejd Hashorva}
\address{Enkelejd Hashorva, Department of Actuarial Science, University of Lausanne, UNIL-Dorigny 1015 Lausanne, Switzerland}
\email{enkelejd.hashorva@unil.ch}

\author{Peng Liu}
\address{Peng Liu, School of Mathematics, Statistics and Actuarial Science, University of Essex, Wivenhoe Park, CO4 3SQ Colchester, UK}
\email{peng.liu@essex.ac.uk}

\bigskip

\date{\today}
 \maketitle
\bigskip

	{\bf Abstract:} We study the asymptotics of sojourn time of the stationary
queueing process $Q(t),t\ge0$ fed by a fractional Brownian motion with Hurst parameter $H\in(0,1)$
above a high threshold $u$.
For the Brownian motion case $H=1/2$, we derive the exact asymptotics of
 \[
 \pk{\int_{T_1}^{T_2}\I(Q(t)>u+h(u))d t>x \Big{|}Q(0) >u }
 \]
 as $u\to\infty$, \eL{where $T_1,T_2, x\geq 0$ and $T_2-T_1>x$}, whereas for all $H\in(0,1)$, we obtain sharp asymptotic approximations of
 \[
 \pk{ \frac 1 {v(u)}  \int_{[T_2(u),T_3(u)]}\I(Q(t)>u+h(u))dt>y \Bigl  \lvert \frac 1 {v(u)} \int_{[0,T_1(u)]}\I(Q(t)>u)dt>x}, \quad x,y >0
 \]
 as $u\to\infty$, for appropriately chosen $T_i$'s and $v$.
Two regimes of the ratio between $u$ and $h(u)$, that lead to qualitatively different approximations, are considered.



\section{Introduction}
Fluid queueing systems with Gaussian-driven structure attained a substantial research interest over the last years; see, e.g., the monograph \cite{Man07} and references therein. Following the seminal contributions \cite{Nor94, Taq97, Mik}
the class of fractional Brownian motions (fBm's) is a
well legitimated model for the traffic stream in modern communication networks.

Let  $B_H(t), t\in\mathbb{R}$ be a standard  fBm with Hurst index $H\in (0,1)$, that is a Gaussian process with
continuous sample paths, zero mean and  covariance function satisfying
$$\eHe{2}Cov(B_H(t), B_H(s))=|s|^{2H}+|t|^{2H}-|t-s|^{2H}, \quad s,t\in\mathbb{R}.$$
Consider the fluid queue fed by $B_H$ and emptied with \eL{a constant rate} $c>0$.
Using the interpretation that for $s<t$, the increment
$B_H(t)-B_H(s)$ models the amount of traffic that entered the system
in the time interval $[s,t)$, we define the workload process $Q(t), t\ge 0$ by
\begin{eqnarray}
Q(t)=B_H(t)-ct+\max\left(Q(0),-\inf_{s\in [0,t]}(B_H(s)-cs)\right).
\label{def.Q}
\end{eqnarray}
The unique stationary solution to the above equation, that is the object of the analysis in this contribution, takes the following form (see. e.g., \cite{Man07})
\BQN\label{queue}
\{Q(t), t\geq 0\}\stackrel{d}{=}\left\{\sup_{s\geq t}(B_H(s)-B_H(t)-c(s-t)), t\geq 0\right\}.
\EQN
\eL{The vast majority} of the analysis of queueing models with Gaussian inputs
deals with the asymptotic results, with particular focus on the
asymptotics of the  probability
\[
\pk{Q(t)>u}
\]
as $u\to\infty$, see e.g., \cite{Nor94,HP99,HP2004,DI2005,Man07,KrzysPeng2015}.
Much less is known about transient characteristics of $Q$, such as
\[
\pk{Q(T)> \omega\Big{|}Q(0)>u},
\]
with a notable exception for the Brownian motion ($H=1/2$).
In particular, in view of  \cite{LiM07}, see also related works \cite{AbW87a,AbW87b, AbW88},
\kdd{it is known that} for $H=1/2$ and $u,\omega,T>0$
\BQN\label{q1}
\pk{Q(T)>\omega\Big{|} Q(0)=u}=\Phi\left(\frac{u-\omega-cT}{\sqrt{T}}\right)+e^{-2c\omega}\Phi\left(\frac{\omega+u-cT}{\sqrt{T}}\right)
\EQN
and
\BQN\label{q2}
\pk{Q(T)>\omega\Big{|} Q(0)>u} &=&-e^{2uc}\Phi\left(\frac{-\omega-u-cT}{\sqrt{T}}\right)+e^{-2c(\omega-u)}\Phi\left(\frac{\omega-u-cT}{\sqrt{T}}\right)\\
&&+ \notag
\Phi\left(\frac{u-\omega-cT}{\sqrt{T}}\right)+e^{-2c\omega}\Phi\left(\frac{\omega+u-cT}{\sqrt{T}}\right),
\EQN
where $\Phi(\cdot)$  denotes the distribution function of a standard Gaussian random variable.
\kdd{Since $Q(0)$ is exponentially distributed for $H=1/2$,
(\ref{q1})-(\ref{q2}) lead to explicit formula for
$\pk{Q(0)>u, Q(T)>\omega}$, which compared with
$\pk{Q(0)>u}\pk{ Q(T)>\omega}$
gives
some insight to the dependence structure of the workload process $Q(t),t\ge0$.}
Since the general case   $H\in(0,1)$ is very complicated,
the findings available in the literature
concern mainly large deviation-type results; see e.g.,
\cite{DeM09} where the asymptotics of
\[
\ln(\pk{Q(0)>pu, Q(Tu)>qu}), \quad u\to \IF
\]
 was derived for $H\in(0,1)$. See also  \cite{AM} for corresponding results   the {\it many-source} model.

\eL{\eHe{In addition to} the \kdd{conditional probability (\ref{q2}),
it is also interesting to know how much time the queue \eHe{spends}  above a \kdd{given}
threshold during a \eHe{given} time period.}}
This motivates us to consider \kdd{the following quantity} 
\begin{eqnarray}
\pk{\int_{T_1}^{T_2}\I(Q(t)>\omega) dt> x \Big{|}Q(0) >u }, \quad   x\in [0, T_2-T_1)
\label{soj.bm}
\end{eqnarray}
for given non-negative $T_1< T_2$.

In Section \ref{s.bm}, for $H=1/2$ we derive exact asymptotics of the above conditional sojourn time by letting
$u,\omega=\omega(u)\to \IF$ in an appropriate way. Specifically, we shall distinguish between two regimes that lead to qualitatively different results:
\begin{enumerate}[(i)]
	\item \kdd{small fluctuation regime:} $\omega=u+w+\eHe{o(1)}, w\in\R$, for which the asymptotics of (\ref{soj.bm}) tends to a positive constant
	as $u\to\infty$;
	\item \kdd{large fluctuation regime:} $\omega=(1+a)u+\eHe{o(u)}, a\in (-1,\infty)$, for which (\ref{soj.bm}) tends to $0$ as $u\to\infty$ with the speed controlled by  $a$
\end{enumerate}
see Propositions \ref{prop.1}, \ref{prop.2} respectively.

Then, in Section \ref{s.fbm} for all $H\in (0,1)$  and $x,y$ non-negative we shall investigate approximations, as $u\to \IF$, of  the following conditional sojourn times probabilities
\BQN
\quad \mathscr{P}^{x,y}_{{T}_1,{T}_2,{T}_3}(\omega,u)
:=\pk{ \frac 1{v(u)} \int_{[T_2(u),T_3(u)]}\I(Q(t)>\omega)dt>y \Bigl  \lvert  \frac 1{v(u)}  \int_{[0,T_1(u)]}\I(Q(t)>u)dt>x},
\label{main}
\EQN
where  $T_i(u), i=1,2,3$ and $\omega=u + h(u),v(u)$ are  suitably chosen  functions,
see assumption ({\bf {T}}). In Theorem \ref{TH1}, complementing the findings of Proposition \ref{prop.1},
we shall determine
\BQN \label{Cxy}
 \limit{u} \mathscr{P}^{x,y}_{{T}_1,{T}_2,{T}_3}(u+au^{2H-1},u)
 \EQN
under some asymptotic restrictions on $T_i(u)$'s and $a\in\R$, which yield a positive and finite limit.  The idea of its proof is based on a modification of recently developed extension of the uniform double-sum technique for
functionals of Gaussian processes \cite{DHLM23}.
Then, in Theorem \ref{TH3} we shall obtain approximations of $\mathscr{P}^{x,y}_{{T}_1,{T}_2,{T}_3}((1+a)u ,u)$ as $u\to\infty$,
which correspond to the results derived in Proposition \ref{prop.2}.
The main findings of this section go in line with recently derived
asymptotics for
\BQNY
\pk{  \frac 1 {v(u)} \int_{[0,T_1(u)]}\I(Q(t)>u)dt>x}
\EQNY
as $u\to\infty$, see \cite{DHLM23}.

{\it Structure of the paper}:
Section \ref{s.bm} is devoted to the analysis of the exact asymptotics of
(\ref{soj.bm}) for the classical model of the Brownian-driven queue, while in Section \ref{s.fbm} we shall investigate asymptotic properties of
$\mathscr{P}^{x,y}_{{T}_1,{T}_2,{T}_3}(\omega ,u)$ for $H\in(0,1)$.
Proofs of all the results are deferred to Section \ref{s.proofs} and Appendix.

\section{Preliminary results}\label{s.bm}
In this section we shall focus on the exact asymptotics
of (\ref{soj.bm}) for the queueing process (\ref{queue}) driven by the Brownian motion.

Let in the following for $T_2-T_1>x\ge 0$ and $w\in\R$
\begin{eqnarray}
\mathcal{C}(T_1,T_2,x;w)=
2c\int^{w}_{-\IF} e^{2cy}\pk{\int_{T_1}^{T_2}\I\left( B_{1/2}(t)-ct >y\right)dt>x}dy\in (0,\IF).\label{C1}
\end{eqnarray}
We begin with a {\it small \kdd{fluctuation}} result concerning the case when the distance between $u$ and \ehe{$\omega=\omega(u)$} in (\ref{soj.bm}) is \ehe{asymptotically} constant. Below the term $o(1)$ is means for $u\to \IF$.
\BS\label{prop.1}
If  $H=1/2$ and $T_2-T_1>x\ge 0$,  then for \ehe{$\omega(u)=u+w+o(1), w\in \R$}
\BQNY
\pk{\int_{T_1}^{T_2}\I(Q(t)>\eHe{\omega(u)})dt>x\Big{|}Q(0)>u}\sim
e^{-2cw}\mathcal{C}(T_1,T_2,x;w)
\EQNY
as $u\to\infty$.
\ES
Next, we consider the {\it large \kdd{fluctuation}} scenario, i.e., $\omega=\omega(u)$ in (\ref{soj.bm}) is asymptotically  proportional to $u$.
\BS\label{prop.2}
Suppose that $H=1/2, T_2-T_1>x\ge 0$ and \eHe{$\omega(u)=  (1+a)u+ o(u)$.}
\begin{enumerate}[(i)]
	\item   If $a\in(-1,0)$, then as $u \to\IF$
\BQNY
\pk{\int_{T_1}^{T_2}\I(Q(t)>\eHe{\omega(u)})dt>x\Big{|}Q(0)>u}
\sim 1.
\EQNY
\item If $a>0$, then as $u \to\IF$
\BQNY
\pk{\int_{T_1}^{T_2}\I(Q(t)>\eHe{\omega(u)})dt>x\Big{|}Q(0)>u}
\sim
e^{-2c(\OU -u)}\mathcal{C}(T_1,T_2,x;\IF).
\EQNY
\end{enumerate}
\ES
 \eL{Both (i) and (ii) in \Cref{prop.2}} also hold if $T_1,T_2$ depend on $u$  \eL{in such a way} that as $u\to \IF$, these converge to positive constants $T_1 < T_2$ with \eL{$T_2-T_1>x\geq 0$}.

\section{Main results}\label{s.fbm}
This section is devoted to the asymptotic analysis of (\ref{main})
for the queueing process $Q$ defined in (\ref{queue}) with \eHe{fBm} input $B_H$, $H\in(0,1)$.
Before proceeding to the main results of this contribution,
we introduce some notation and assumptions.
 Let $W_H(t)= \sqrt{2}B_H(t)-|t|^{2H}, t\in \R$ and
define for
$$x\ge0 ,y\ge 0$$
 and $\eeh{\lambda\inr}, \mathscr{T}_1>0$, $0<\mathscr{T}_2<\mathscr{T}_3< \IF $
$$ \overline {\mathcal{B}}_H^{x,y}(\mathscr{T}_1;\eeh{\lambda},\mathscr{T}_2,\mathscr{T}_3)= \int_{\R} e^z \pk{\int_{[0,\mathscr{T}_1]}\I( W_H(t)>z)dt>x, \int_{[\mathscr{T}_2,\mathscr{T}_3]}\I( W_H(t)>z+\eeh{\lambda})dt>y}dz $$
and set
$$ \overline{\mathcal{B}}_H^{x}(\mathscr{T}_1) = \int_{\R} e^z \pk{\int_{[0,\mathscr{T}_1]}\I( W_H(t)>z)dt>x}dz.$$
Further, given $H\in (0,1)$, $c>0, u>0$ let
\BQN\label{AB}
A=\left(\frac{{H}}{c(1-{H})}\right)^{-H}\frac{1}{1-H},
\quad t^*=\frac{H}{c(1-H)}, \quad \Delta(u)=2^{\frac{1}{2H}}t^*A^{-\frac{1}{H}}u^{-\frac{1-H}{H}}
\EQN
and set
\[
v(u)=u\Delta(u).
\]
In the rest of this section, for a given function $h$, we analyse the asymptotics of \eL{ $\mathscr{P}^{x,y}_{{T}_1,{T}_2,{T}_3}(\eHe{\omega(u)},u)$ defined in \eqref{main} with $\eHe{\omega(u)}=u+h(u)$}
as $u\to\infty$,
where  \eL{$T_i(u)$'s depend on  $u$ in such a way} that \\
\\
{\bf (T)} ${\limit{u}}\frac{T_i(u)}{v(u)}=\mathscr{T}_i\in (0,\IF)$, for $i=1,2,3$ with
$\mathscr{T}_1>x$ and  $\mathscr{T}_3-\mathscr{T}_2>y$\\
\\
is satisfied.

\eL{We note in passing that for $H=1/2$, $v(u)=u\Delta(u)=2^{\frac{1}{2H}}t^*A^{-\frac{1}{H}}$ is a constant}. Hence,  under  {\bf (T)}, we have $T_i(u)\to \mathscr{C}_i\in(0,\infty)$. Thus {\bf (T)}
included the model considered in Propositions \ref{prop.1} and \ref{prop.2}.

We shall consider two scenarios that depend on the relative size of $h(u)$ with respect to $u$:
\begin{itemize}
  \item[$\diamond$] {\it small \kdd{fluctuation}} case:  $|h(u)|$ is relatively {\it small} with respect to $u$, i.e., \eL{$h(u)=\lambda u^{2H-1}$ with $\lambda\in \R$ and $H\in (0,1)$},  which leads to
$\lim_{u\to\infty} \mathscr{P}^{x,y}_{{T}_1,{T}_2,{T}_3}(u+h(u),u) > 0$;
  \item[$\diamond$] {\it large \kdd{fluctuation}} case:  $h(u)=au$ is \dd{proportional to}
$u$, \dd{which leads to}    $\mathscr{P}^{x,y}_{{T}_1,{T}_2,{T}_3}(u+h(u),u)\to 0$ if $h(u)>0$ \dd{and}
$\mathscr{P}^{x,y}_{{T}_1,{T}_2,{T}_3}(u+h(u),u)\to 1$ if $h(u)<0$ as $u\to\infty$.
\end{itemize}

{\it Small \kdd{fluctuation} regime.}
We begin with the case when $h(u)$ is relatively small with comparison
to $u$ and thus the conditional probability
$\mathscr{P}^{x,y}_{{T}_1,{T}_2,{T}_3}(u+h(u),u)$
is cut away from $0$, as $u\to\infty$.

\BT\label{TH1} If {\bf (T)} holds, then with $Q$ defined in (\ref{queue}) and $\lambda\inr$
\BQN
\lim_{u\to\infty}
\mathscr{P}^{x,y}_{{T}_1,{T}_2,{T}_3}\left(u+ {\frac{\lambda}{\dd{A^2}(1-H)} u^{2H-1}},u\right)=  \frac{ \overline{\mathcal{B}}_H^{x,y}(\mathscr{T}_1;\eeh{\lambda}, \mathscr{T}_2,\mathscr{T}_3)}{\overline{\mathcal{B}}_H^{x}(\mathscr{T}_1)}
\in (0,\IF).\label{l1}
\EQN
\ET

\begin{remark}
	\begin{enumerate}[(i)]
		\item \eeh{In the case of Brownian motion with  $H=1/2$, function
$v(u)= 1/(2 c^2)$ does not depend on $u$ and the above reads}
		\BQN
		\lim_{u\to\infty}
		\mathscr{P}^{x,y}_{{T}_1,{T}_2,{T}_3}\left(u+ \rd{\frac{\lambda}{{\dd{2}}c}},u\right)
		&=&
		\frac{ \overline{\mathcal{B}}_H^{x,y}(\mathscr{T}_1;\lambda, \mathscr{T}_2,\mathscr{T}_3)}
{\overline{\mathcal{B}}_H^{x}(\mathscr{T}_1)} 		\in (0,\IF).\label{l.1}
		\EQN
\eeh{		
		Since $v(u)$ is constant in this case, we can take $T_i = 2c^2\mathscr{T}_i>0,i\le 3$ \dd{in (\ref{l.1})}.
In the particular case that $x=0, H=1/2$ we have
		\BQNY
\limit{u} 		\pk{\int_{[T_2,T_3]}\I\left(Q(t)>u+{\frac{\lambda}{{\dd{2}}c}}\right)dt> 2 c^2 y \Bigl  \lvert \sup_{t\in [0,T_1]}Q(t)>u }
		&=& 		\frac{ \overline{\mathcal{B}}_H^{0,y}(\mathscr{T}_1;\lambda, \mathscr{T}_2,\mathscr{T}_3)}{\overline{\mathcal{B}}_H^{0}(\mathscr{T}_1)} 		\in (0,\IF)
		\EQNY
and taking $y=0$ yields
		\BQNY
\limit{u} 		\pk{\sup_{t\in [T_2,T_3]}Q(t)>u+\frac{\lambda}{{\dd{2}}c}  \Bigl  \lvert \sup_{t\in [0,T_1]}Q(t)>u }
&=& 		\frac{ \overline{\mathcal{B}}_H^{0,0}(\mathscr{T}_1;\lambda, \mathscr{T}_2,\mathscr{T}_3)}{\overline{\mathcal{B}}_H^{0}(\mathscr{T}_1)} 		\in (0,\IF).
\EQNY
}		

		\item It follows from \Cref{TH1} that
		for $h(u)=o(u^{2H-1})$
		\BQN
		\lim_{u\to\infty}
		\mathscr{P}^{x,y}_{{T}_1,{T}_2,{T}_3}(u+h(u),u)=
        \frac{\dd{ \overline{\mathcal{B}}_H^{x,y}(\mathscr{T}_1; 0 , \mathscr{T}_2,\mathscr{T}_3)}}{\overline{\mathcal{B}}_H^{x}(\mathscr{T}_1)} \in (0,\IF).
		\label{lrd3}
		\EQN
\dd{		Notably, if $H \in (1/2,1)$, then
$T_i(u)\sim  \mathscr{T}_i u^{(2H-1)/H}$ as $u\to\infty$ for $i=1,2,3$.
Hence
$\lim_{u\to\infty}(T_2(u)-T_1(u))=\infty$
and one can take $h(u)\to\infty$, as $u\to\infty$. Thus, the insensitivity of
limit (\ref{lrd3}) on $h(u)$  is yet another manifestation of the
{\it long range dependence} property of  $Q$  inherited from
the input process $B_H$.
This observation  goes in line with the {\it{Piterbarg property}}
\[
\lim_{u\to\infty} \frac{\pk{\sup_{t\in [0,T(u)]} Q(t)>u}}{\pk{Q(0)>u}}=1
\]
derived in \cite{Pit01} and
{\it the strong Piterbarg property} see \cite{DeK14}, namely
\[
\lim_{u\to\infty} \frac{\pk{\inf_{t\in [0,T(u)]} Q(t)>u}}{\pk{Q(0)>u}}=1,
\]
where $T(u)=o(u^{(2H-1)/H})$ as $u\to\infty$.
 }
	\end{enumerate}
\end{remark}

{\it Large \kdd{fluctuation} regime.}
Suppose next that $h(u)=a u$, $a\neq0$.
It appears that in this case the \kdd{fluctuation} $h(u)$ substantially influences the asymptotics of
$\mathscr{P}^{x,y}_{{T}_1,{T}_2,{T}_3}(u+h(u),u)$ as $u\to\infty$.
We point out the lack of symmetry with respect to the sign of $a$ in the results given in the following theorem,
\dd{which is due to the non-reversibility in time of the queueing process $Q$},
\kdd{i.e., the fact that
	$$\pk{Q(s)>u,Q(t)>v}\neq \pk{Q(t)>u,Q(s)>v}$$ for $u\neq v$}.
\BT\label{TH3} Let $Q$ be defined in (\ref{queue}) and \eHe{set $\tilde{a}=(1+a)^{(1-2H)/H}$}. Suppose that   {\bf (T)} holds.
\begin{enumerate}[(i)]
 \item If $a\in(-1,0)$, then
\BQNY
\lim_{u\to\infty}\mathscr{P}^{x,y}_{{T}_1,{T}_2,{T}_3}((1+a)u,u)=1.
\EQNY
\item  If $a>0$, then
\BQNY
\limsup_{u\to\infty}
\frac{\mathscr{P}^{x,y}_{{T}_1,{T}_2,{T}_3}((1+a)u,u)}
{\exp\left(-\frac{A^2\left((1+a)^{2-2H}-1\right)}{2}u^{2-2H}\right)}
\le
\tilde{a}^{1-H}\frac{\overline{\mathcal{B}}_H^{\tilde{a}y}(\tilde{a}(\mathscr{T}_3-\mathscr{T}_2))}
{\overline{\mathcal{B}}_H^{x}(\mathscr{T}_1)}
\EQNY
and
\BQNY
\liminf_{u\to\infty}
\frac{\mathscr{P}^{x,y}_{{T}_1,{T}_2,{T}_3}((1+a)u,u)}
{\exp\left(-\frac{A^2\left((1+a)^{2-2H}-1\right)}{2}u^{2-2H}\right)}
\ge
\tilde{a}^{1-H}
\frac{
\overline{\mathcal{B}}_H^{\tilde{a}x,\tilde{a}y}(\tilde{a}\mathscr{T}_1;\eeh{0},\tilde{a} \mathscr{T}_2,\tilde{a}\mathscr{T}_3)}
{\overline{\mathcal{B}}_H^{x}(\mathscr{T}_1)}.
\EQNY
\end{enumerate}
\ET
\dd{
\begin{remark}
Theorem \ref{TH3} straightforwardly implies that
\[
\lim_{u\to\infty}\frac{\ln\left(  \mathscr{P}^{x,y}_{{T}_1,{T}_2,{T}_3}((1+a)u,u)  \right)}{u^{2-2H}}
=
-\frac{1}2 A^2\left((1+a)^{2-2H}-1\right), \quad \forall a>0.
\]
\end{remark}
}

\section{Proofs}\label{s.proofs}
In this section we present detailed proofs of \Cref{prop.1}, \ref{prop.2}
and \Cref{TH1}, \ref{TH3}.

\subsection{Proof of Proposition \ref{prop.1}}
Recall that
by (\ref{def.Q})
\[
Q(t)=B_{1/2}(t)-ct+\max\left(Q(0),-\inf_{s\in [0,t]}(B_{1/2}(s)-cs)\right),
\]
where
$Q(0)$ is independent of $B_{1/2}(t)-ct$ and $\inf_{s\in [0,t]}(B_{1/2}(s)-cs)$ for $t> 0$.
\eHe{By   \cite[Eq.\ (5)]{DLM18}  we have}
\BQN
\pk{Q(0)>u}=\pk{\sup_{t\geq 0} (B_{1/2}(t)-ct)>u}=e^{-2cu}, \quad u\geq 0.
\label{exp}
\EQN
Hence it suffices to analyse
\BQNY
\pk{\int_{T_1}^{T_2}\I(Q(t)>\OU)dt>x, Q(0)>u}.
\EQNY
We note first that
\BQN
\lefteqn{\pk{\int_{T_1}^{T_2}\I(Q(t)> \OU )dt>x, Q(0)>u}}\nonumber\\
&&\ge
\pk{\int_{T_1}^{T_2}\I(Q(0)+B_{1/2}(t)-ct>\OU )dt>x, Q(0)>u}.\label{low1}
\EQN
Moreover, we have
\BQN
\lefteqn{\pk{\int_{T_1}^{T_2}\I(Q(t)>\OU)dt>x, Q(0)>u}}\nonumber\\
&=&
\pk{\int_{T_1}^{T_2}\I(Q(t)>\OU)dt>x, Q(0)>u, \rd{\sup_{s\in [0,T_2]}}(cs- B_{1/2}(s)) \le u}\nonumber\\
&&+
\pk{\int_{T_1}^{T_2}\I(Q(t)>\OU)dt>x, Q(0)>u, \rd{\sup_{s\in [0,T_2]}}(cs- B_{1/2}(s)) > u}\nonumber\\
&=&
P_1(u)+P_2(u).\label{upp1}
\EQN
For $P_1(u)$  we have the following upper bound
\begin{eqnarray*}
P_1(u)\le
\pk{\int_{T_1}^{T_2}\I(Q(0)+B_{1/2}(t)-ct>\OU)dt>x, Q(0)>u}
\end{eqnarray*}
and for $P_2(u)$  by Borell-TIS inequality (see, e.g., \cite{AdlerTaylor})
\begin{eqnarray}
P_2(u)
&\le&
\pk{\sup_{s\in [0,T_2]}(cs- B_{1/2}(s)) >u}
\le
e^{-Cu^2}\label{square}
\end{eqnarray}
for some $C>0$ and sufficiently large $u$.

Next, we note that
\begin{eqnarray*}
\lefteqn{
\pk{\int_{T_1}^{T_2}\I(Q(0)+B_{1/2}(t)-ct>\OU)dt>x, Q(0)>u}}\\
&&=
2c\int_{u}^\IF e^{-2cy}\pk{\int_{T_1}^{T_2}\I\left(y+B_{1/2}(t)-ct>\OU \right)dt>x}dy\\
&&=
2ce^{-2c\OU}
\int_{u -\OU}^\IF e^{-2cy}\pk{\int_{T_1}^{T_2}\I\left( B_{1/2}(t)-ct >-y\right)dt>x}dy\\
&&=
2ce^{-2c\OU}
\int^{\OU-u}_{-\IF} e^{2cy}\pk{\int_{T_1}^{T_2}\I\left( B_{1/2}(t)-ct >y\right)dt>x}dy
\\
&&=
e^{-2c\OU }\mathcal{C}(T_1,T_2,x;\OU-u),
\end{eqnarray*}
where $\mathcal{C}(T_1,T_2,x;z)$ is defined in (\ref{C1}).
Hence, by (\ref{square}) applied to (\ref{low1}) and (\ref{upp1}),
we arrive at
\begin{eqnarray}
\pk{\int_{T_1}^{T_2}\I(Q(t)>\OU)dt>x, Q(0)>u}
\sim
e^{-2cu}e^{-2cw}\mathcal{C}(T_1,T_2,x;w)
\label{sim.BM}
\end{eqnarray}
as $u\to\infty$.
Finally, by (\ref{exp}) we get
\[
\pk{\int_{T_1}^{T_2}\I(Q(t)>\OU)dt>x \big{|} Q(0)>u}\sim
e^{-2cw}\mathcal{C}(T_1,T_2,x;w)
\]
as $u\to\infty$. This completes the proof.
\QED

\subsection{Proof of Proposition \ref{prop.2}}
The idea of the proof is the same as the proof of Proposition \ref{prop.1}.
\eHe{Since} by Borell-TIS inequality
\begin{eqnarray*}
\pk{\int_{T_1}^{T_2}\I\left( B_{1/2}(t)-ct >y\right)dt>x}
&\le&
\pk{\sup_{t\in [T_1,T_2]}\left( B_{1/2}(t)-ct\right)>y}\\
&\le&
C_1 \exp(-C_2y^2)
\end{eqnarray*}
for some positive constants $C_1,C_2$, we conclude that
\[
\mathcal{C}(T_1,T_2,x;\IF)
=
\int^{\IF}_{-\IF} e^{2cy}\pk{\int_{T_1}^{T_2}\I\left( B_{1/2}(t)-ct >y\right)dt>x}dy<\infty.
\]
Thus
if $a\in(-1,0)$,  then as $u\to \IF$
\BQNY
&&\pk{\int_{T_1}^{T_2}\I(Q(t)>\OU )dt>x, Q(0)>u}\\
&&\quad \sim
2ce^{-2c\OU}
\int^{\OU- u}_{-\IF} e^{2cy}\pk{\int_{T_1}^{T_2}\I\left( B_{1/2}(t)-ct >y\right)dt>x}dy\\
&&\quad \sim
e^{-2cu},
\EQNY
where we used that uniformly for $y\in (-\infty, \OU- u]$
\[
\limit{u} \pk{\int_{T_1}^{T_2}\I\left( B_{1/2}(t)-ct >y\right)dt>x}=1.
\]
\eHe{Similarly} for
$a>0$, we have that
\BQNY
\pk{\int_{T_1}^{T_2}\I(Q(t)>\OU)dt>x, Q(0)>u}
 \sim e^{-2c\OU}
\mathcal{C}(T_1,T_2,x;\IF)
\EQNY
as $u\to\infty$. Thus, combining the above with (\ref{exp}), we complete the proof.
\QED

\subsection{Proof of Theorem \ref{TH1}}
\eeh{We begin with a result which is crucial for the proof and of some interests on its own right.}
Recall that $Q$ is defined in (\ref{queue}).
For $B_H, B_H'$ two independent fBm's with Hurst indexes $H$, we set
\begin{eqnarray}
	W_H(t)= \sqrt{2}B_H(t)-|t|^{2H}, \quad W_H'(t)= \sqrt{2}B_H'(t)-|t|^{2H}
	\label{eWH}
\end{eqnarray}
and
$$
V_H(S)=  \sup_{s\in [0,S]} W_H(s), \quad t\inr, S\ge 0.$$
Define further \Eh{for all $x,y$ non-negative} and $\lambda\in\R$, the generalized Berman-type constants  by
\BQNY &&\mathcal{B}_H^{x,y}(T_1; \lambda, T_2,T_3)([0,S])\\
&&=\int_{\R} e^w \pk{\int_{[0,T_1]}\I( W_H'(t)+
	V_H(S)>w)dt>x,
	\int_{[T_2,T_3]}\I( W'_H(t)+  V_H(S)
>w+\lambda)dt>y}dw.
\EQNY
Denote further by $ \mathcal{H}_{2H}$ the Pickands constant corresponding to $B_H$, i.e.,
$$  \mathcal{H}_{2H} =\eL{\lim_{S\to\infty}S^{-1} {\E{e^{V_H(S)} }}}=\E{ \frac{ \sup_{t\in \R} e^{ W_H(t)} } { \int_{t\in \R} e^{ W_H(t)}  dt}} \in (0,\IF).$$
\BEL\label{lemA0} For all $T_1,T_2,T_3$ positive, $\lambda\in\R$, and all $x,y$ non-negative we have
\BQNY
\mathcal{B}_H^{x,y}(T_1;\lambda,T_2,T_3):=\lim_{S\rw\IF}S^{-1} \mathcal{B}_H^{x,y}(T_1;\lambda,T_2,T_3)([0,S])
&=&  \mathcal{H}_{2H}  \overline {\mathcal{B}}_H^{x,y}(T_1;\lambda, T_2,T_3) \in (0,\IF).
\EQNY
\EEL
\eL{It is worth mentioning that both sides of equation in the above lemma is equal to zero if $x\geq T_1$ or $y\geq T_3-T_2$. Hence it is valid for all nonnegative $x$ and $y$.}

\prooflem{lemA0} First note that for any $S>0$ we have using the Fubini-Tonelli
theorem and the independence of $V_H$ and $W_H'$
\BQNY
\lefteqn{\mathcal{B}_H^{x,y}(T_1;\lambda,T_2,T_3)([0,S])}\\
 &=&  \E{ \int_\R e^w
	\I( \int_{[0,T_1]}\I( W_H'(t)+ V_H(S)>w)dt>x, \int_{[T_2,T_3]}\I( W'_H(t)+  V_H(S)>w+\lambda)dt>y )dw}\\
&=&\E{e^{V_H(S)} \int_\R e^{w}
	\I( \int_{[0,T_1]}\I( W_H'(t)>w)dt>x, \int_{[T_2,T_3]}\I( W'_H(t)>w+\lambda)dt>y )dw}\\
&=&\E{e^{V_H(S)} } \int_\R e^{w}
\pk{  \int_{[0,T_1]}\I( W_H(t)>w)dt>x, \int_{[T_2,T_3]}\I( W_H(t)>w+\lambda)dt>y }dw\\
&\le &\E{e^{V_H(S)} } \int_\R e^{w}
\pk{  \int_{[0,T_1]}\I( W_H(t)>w)dt>0 }dw\\
&= &\E{e^{V_H(S)} } \int_0^\IF  e^{w}
\pk{ V_H(T_1)>w }dw\\
&= &\E{e^{V_H(S)} } \E{e^{V_H(T_1)} }.
\EQNY
Hence the  claim follows by the definition of the Pickands constant and the sample continuity of $V_H$.
\QED
\def\ru{u}
\def\uA{A \ru^{1-H}}
Let in the following
$$B=\left(\frac{H}{c(1-H)}\right)^{-H-2}H$$
\eeh{and recall that
	\BQN \label{recD} \Delta(u)=2^{\frac{1}{2H}}t^*A^{-\frac{1}{H}}u^{-\frac{1-H}{H}}, \quad
	v(u)= u \Delta(u).
\EQN
}
Applying \cite[Lem 4.1]{DHLM23} we obtain the following result.
\BS\label{asymptotics}
If {\bf (T)} holds, then
\BQN\label{p1}
\pk{  \frac 1 {v(u)}\int_{[0,T_1(u)]}\I(Q(t)>u)dt>x}\sim \mathcal{H}_{2H}\overline{\mathcal{B}}_H^{x}(\mathscr{T}_1)\frac{\sqrt{2\pi}(AB)^{-1/2}}
{u^{1-H}\Delta(u)}\Psi(\uA),\quad u\rw\IF.
\EQN
\ES
The next proposition   \eL{plays a key role} in the proof of Theorem \ref{TH1}.
\def\ruu{u+ \tau u^{2H-1}}
\BS\label{double.as}
If {\bf (T)} holds, then for all $\lambda\in \R, \dd{\tau=\lambda/(A^2(1-H))}$ 
\BQN\label{p2}
\lefteqn{\pk{ \nvu  \int_{[T_2(u),T_3(u)]}\I(Q(t)> \ruu )dt>y, \nvu  \int_{[0,T_1(u)]}\I(Q(t)>\ru)dt>x}}\nonumber\\
&&\sim \mathcal{H}_{2H}\overline{\mathcal{B}}_H^{x,y}(\mathscr{T}_1;\eeh{\lambda},\mathscr{T}_2,\mathscr{T}_3)\frac{\sqrt{2\pi}(AB)^{-1/2}}
{\rd{u^{1-H}}\Delta(u)}\Psi(\uA),\quad u\rw\IF.
\EQN
\ES

\def\ruuu{ \eeh{\widetilde{u}}}
\def\uA{{u_\star}}
\def\uAA{{\widetilde{u}_\star}}

Hereafter, for any \eL{non-constant random variable} $Z$, we denote $\overline{Z}=Z/\sqrt{Var(Z)}$.

\proofprop{double.as}
Using the self-similarity of $B_H$, i.e.,
$$\{B_H(ut), t\in \mathbb{R}\}\stackrel{d}{=}\{u^HB_H(t), t\in\mathbb{R}\},\quad u>0
$$
 we have with $\Delta(u)$ given in \eqref{AB} and
$\ruuu= \ruu$
\BQNY
&&\pk{\nvu \int_{[0,T_1(u)]}\I(Q(t)>\ru)dt>x, \nvu \int_{[T_2(u),T_3(u)]}\I(Q(t)>\ruuu)dt>y}\\
&&\quad=\pk{\ndu \int_{[0,T_1(u)/\ru]}\I(\sup_{s\geq t}(\ru^H(B_H(s)-B_H(t))-c\ru(s-t))>\ru)dt>x,\right.\\
&&\quad \quad \left. \ndu \int_{[T_2(u)/\dd{\ruuu},T_3(u)/\dd{\ruuu}]}\I(\sup_{s\geq t}(\dd{\ruuu}^H(B_H(s)-B_H(t))-c\ruuu(s-t))>\ruuu)dt>y}\\
&&\quad=\pk{\ndu \int_{[0,\tnu{1} ]}\I(\sup_{s\geq t}Z(s,t)>\uA)dt>x,\right.
\\
&&\quad\quad \left. \ndu  \int_{[\tnu{2},\tnu{3}]}\I(\sup_{s\geq t}Z(s,t)> \uAA )dt>y},
\EQNY
where
$$Z(s,t)=A\frac{B_H(s)-B_H(t)}{1+c(s-t)}$$
and
$$\uA= A \ru^{1-H},  \quad \uAA=A \dd{\ruuu^{1-H}}, \quad \dd{\tnu{1}= T_{1}(u)/\ru}, \dd{\tnu{i}= T_{i}(u)/\ruuu}, i=2,3.$$
Note that as $u\to \IF$
\begin{eqnarray}
	\label{ciediamo}\uAA
	=  \dd{\uA + \frac{\lambda}{\uA}}, \quad \uAA^2\sim \uA^2+ 2 \lambda + o(1) .
\end{eqnarray}
Direct calculation shows that
\[\max_{s\ge t}\sqrt{ Var(Z(s,t))} = \max_{s\ge t} \frac{A(s-t)^{H}}{1+c(s-t)}   = 1\]
and the maximum is attained for all $s,t$ such that
$$s-t=t^*=\frac{H}{c(1-H)}$$
 and
\BQN\label{Var}
1-A\frac{t^H}{1+ct}\sim\frac{B}{2A}(t-t^*)^2, \quad t\rw t^*.
\EQN
Moreover, we have
\BQN\label{Cor}
\lim_{\delta\rw 0} \sup_{|s-t-t^*|, |s'-t'-t^*|<\delta, |s-s'|<
\delta}\left|\frac{1-Cor(Z(s,t), Z(s',t'))}{|s-s'|^{2H}+|t-t'|^{2H}}-2^{-1}(t^*)^{-2H}\right|=0.
\EQN
In the following we tacitly assume that $$S>\max(x,y).$$
Observe that
\BQNY
\pi_1(u)
&\leq&
\pk{\ndu\int_{[0,\tnu{1}]}\I(\sup_{s\geq t}Z(s,t)>\uA)dt>x,  \ndu\int_{[\tnu{2},\tnu{3}]}\I(\sup_{s\geq t}Z(s,t)>\uAA)dt>y}\\
  &\leq& \pi_1(u)+\pi_2(u),
\EQNY
where
\BQNY\pi_1(u)&=&\pk{\ndu \int_{[0,\tnu{1}]}\I(\sup_{|s-t^*|\leq (\ln u)/u^{1-H}}Z(s,t)>\uA)dt>x, \right.\quad\\
&&\quad\quad  \left. \ndu \int_{[\tnu{2},\tnu{3}]}\I(\sup_{|s-t^*|\leq (\ln u)/u^{1-H}}Z(s,t)>\uAA)dt>y},
\EQNY
\def\tnuu{\eeh{\overline{T^*}(u)}}
$$\pi_2(u)=\pk{\sup_{t\in [0,\tnuu]}\sup_{|s-t^*|\geq (\ln u)/(2u^{1-H}), s\geq t}Z(s,t)>\hat{u}},$$
with $\tnuu=\max(\tnu{1}, \tnu{2}, \tnu{3})$ and $\hat{u}=\min(\uA,\uAA)$.\\
{\it$\diamond$ \underline{Upper bound of $\pi_2(u)$}}.
Next, for some $T>0$ we have
\BQNY
\pi_2(u)\leq\pk{\sup_{t\in [0,\tnuu]}\sup_{|s-t^*|\geq (\ln u)/(2u^{1-H}), t\leq s\leq T}Z(s,t)>\hat{u}}+\pk{\sup_{t\in [0,\tnuu]}\sup_{ s\geq T}Z(s,t)>\hat{u}}.
\EQNY
In view of (\ref{Var}) for $u$ sufficiently large
$$
\sup_{t\in [0,\tnuu]}\sup_{|s-t^*|\geq (\ln u)/(2u^{1-H}), t\leq s\leq T}Var(Z(s,t))\leq 1-\mathbb{Q}\left(\frac{\ln u}{u^{1-H}}\right)^2
$$
and by (\ref{Cor})
\BQNY
\E{\left(Z(s,t)-Z(s',t')\right)^2}\leq \mathbb{Q}_1(|s-s'|^H+|t-t'|^H), \quad t\in [0,\tnuu], |s-t^*|\geq (\ln u)/(2u^{1-H}), t\leq s\leq T.
\EQNY
Hence, in light of \cite[Thm 8.1]{Pit96} for all $u$ large enough
\BQNY
\pk{\sup_{t\in [0,\tnuu]}\sup_{|s-t^*|\geq (\ln u)/(2u^{1-H}), t\leq s\leq T}Z(s,t)>\hat{u}}\leq \mathbb{Q}_2u^{\frac{4(1-H)}{H}}\Psi\left(\frac{\hat{u}}{\sqrt{1-\mathbb{Q}\left(\frac{\ln u}{u^{1-H}}\right)^2}}\right).
\EQNY
Moreover, for $T$ sufficiently large
\BQNY
\sqrt{Var(Z(s,t))}=\frac{A(s-t)^H}{1+c(s-t)}\leq \frac{2A}{c}(T+k)^{-(1-H)}, \quad
s\in [T+k, T+k+1], t\in [0,\tnuu].
\EQNY
Hence for some \eeh{$\ve\in (0,1)$ \eHe{(set $c_\ve=(1+\ve)c)$}}
\BQNY
\pk{\sup_{t\in [0,\tnuu ]}\sup_{ s\geq T}Z(s,t)>\hat{u}}&\leq& \sum_{k=0}^\IF\pk{\sup_{t\in [0,\tnuu]}\sup_{ s\in[T+k, T+k+1]}Z(s,t)>\hat{u}}\\
&\leq & \sum_{k=0}^\IF\pk{\sup_{t\in [0,\tnuu ]}\sup_{ s\in[T+k, T+k+1]}\overline{Z}(s,t)>\frac{1}{2}c_\ve(T+k)^{(1-H)}u^{1-H}}.
\EQNY
Additionally, for $T$ sufficiently large and $k\geq 0$,
we have
\BQNY
\E{(\overline{Z}(s,t)-\overline{Z}(s',t'))^2}\leq \mathbb{Q}_3(|s-s'|^{H}+|t-t'|^H), \quad s, s'\in[T+k, T+k+1], t, t'\in [0,1].
\EQNY
Thus by \cite[ Thm 8.1]{Pit96} for all $T$ and $u$ sufficiently large we have
\BQNY
\pk{\sup_{t\in [0,\dd{\overline{T^*}(u)}]}\sup_{ s\geq T}Z(s,t)>\hat{u}}
&\leq & \sum_{k=0}^\IF\pk{\sup_{t\in [0,1]}\sup_{ s\in[T+k, T+k+1]}\overline{Z}(s,t)>\frac{1}{2}c_\ve(T+k)^{(1-H)}u^{1-H}}\\
&\leq&\sum_{k=0}^\IF \mathbb{Q}_4u^\frac{4(1-H)}{H} \Psi\left(\frac{1}{2}c_\ve(T+k)^{(1-H)}u^{1-H}\right)\\
&\leq& \sum_{k=0}^\IF \mathbb{Q}_4u^\frac{4(1-H)}{H} e^{-\frac{1}2 \left(\frac{1}{2}c_\ve(T+k)^{(1-H)}u^{1-H}\right)^2	}\\
&\leq& \sum_{k=0}^\IF \mathbb{Q}_4u^\frac{4(1-H)}{H}\int_{T-1}^{\infty} e^{-\frac{1}2 \left(\frac{1}{2}c_\ve z^{(1-H)}u^{1-H}\right)^2}dz\\
&\leq & \mathbb{Q}_4u^\frac{4(1-H)}{H} \Psi\left(\mathbb{Q}_5(Tu)^{1-H}\right).
\EQNY
Therefore we conclude that for all $u,T$  sufficiently large
\BQN\label{pi2}
\pi_2(u)\leq \mathbb{Q}_2u^{\frac{4(1-H)}{H}}\Psi\left(\frac{\hat{u}}{\sqrt{1-\mathbb{Q}\left(\frac{\ln u}{u^{1-H}}\right)^2}}\right)+ \mathbb{Q}_4u^\frac{4(1-H)}{H} \Psi\left(2Au ^{1-H}\right).
\EQN
{\it $\diamond$ \underline{Upper bound of $\pi_1(u)$}}.
Given a positive integer $k$ and $u>0$ define
$$\quad I_{k}(u)=[k\Delta(u)S, (k+1)\Delta(u)S], \quad  N(u)=\left[\frac{\ln u}{u^{1-H}\Delta(u)S}\right]+1.
$$
It follows that
\BQNY
\pi_1(u)&=&\pk{\ndu \int_{[0,\tnu{1}]}\I(\sup_{|s|\leq (\ln u)/u^{1-H}}Z(s+t^*,t)>\uA)dt>x,\right.\\
&&\quad \left. \ndu \int_{[\tnu{2},\tnu{3}]}\I(\sup_{|s|\leq (\ln u)/u^{1-H}}Z(s+t^*,t)>\uAA)dt>y}\\
&\leq& \Sigma_1^+(u)+2\Sigma\Sigma_1(u)+2\Sigma\Sigma_2(u),
\EQNY
where
\BQNY
\Sigma_1^+(u)&=&\sum_{k=-N(u)-1}^{N(u)+1}
\pk{\ndu \int_{[0,\tnu{1}]}\I(\sup_{s\in I_k(u)}Z(s+t^*,t)>\uA)dt>x,\right.\\
&&\left.
 \ndu \int_{[\tnu{2},\tnu{3}]}\I(\sup_{s\in I_k(u)}Z(s+t^*,t)>\uAA)dt>y}\\
 &\leq &\sum_{k=-N(u)-1}^{N(u)+1}
 \pk{\int_{[0,\mathscr{T}_1+\epsilon]}\I(\sup_{s\in [0,S]}Z_{u,k}(s,t)>u_k^-)dt>x,\right.\\
&&\left.
 \int_{[\mathscr{T}_2-\epsilon, \mathscr{T}_3+\epsilon]}\I(\sup_{s\in [0,S]}Z_{u,k}(s,t)>\widetilde{u}_k^-)dt>y},\\
 \Sigma\Sigma_1(u)&=&\sum_{|k|,|l|\leq N(u)+1, l=k+1}
 \pk{\sup_{t\in [0,T^*], s\in [kS,(k+1)S]}Z(\Delta(u)s+t^*,\Delta(u)t)>\uAA,\right.\\
 && \left.\sup_{t\in [0, T^*], s\in [lS, (l+1)S]}Z(\Delta(u)s+t^*,\Delta(u)t)>\uA},\\
 \Sigma\Sigma_2(u)&=&\sum_{|k|,|l|\leq N(u)+1, l\geq k+2}\pk{\sup_{t\in [0,T^*], s\in [kS,(k+1)S]}Z(\Delta(u)s+t^*,\Delta(u)t)>\uAA,\right.\\
 && \left.\sup_{t\in [0,T^*], s\in [lS, (l+1)S]}Z(\Delta(u)s+t^*,\Delta(u)t)>\uA},
\EQNY
with
$$ T^*=\max(\mathscr{T}_1+\epsilon,\mathscr{T}_2-\epsilon,\mathscr{T}_3+\epsilon), \epsilon<\mathscr{T}_2, \quad
\Delta(u)=C u^{-\frac{1-H}{H}}, \quad C=2^{\frac{1}{2H}}t^*A^{-\frac{1}{H}},$$
$$Z_{u,k}(s,t)=\overline{Z}(t^*+\Delta(u)(kS+s),\Delta(u)t),$$
$$\quad u_{k}^-=\uA\left(1+  \frac{(1- \epsilon)B}{2A}\Delta^2(u)\eeh{\eta_{k,S}}\right), \quad \eeh{\eta_{k,S}}=\inf_{s\in [kS,(k+1)S], t\in [0, T_*]}(s-t)^2,$$
$$\quad \widetilde{u_{k}}^-=\uAA \left(1+\frac{(1- \epsilon)B}{2A}\Delta^2(u)\eeh{\eta_{k,S}}\right).$$
Since the maximal value of $k$ is $N(u)=\left[\frac{\ln u}{u^{1-H}\Delta(u)S}\right]+1$ and
$\eta_{k,S}$ is non-negative and bounded up to some constant by $k^2S^2$ using further \eqref{ciediamo} we have
\begin{eqnarray}\label{caj}
	u_{k}^-=\uA(1+ o( u^{H-1}\ln u )) , \quad
\widetilde{u_{k}}^-=(\uA+  \lambda/\uA) ( 1+ o( u^{H-1}\ln u ))= \eeh{u_{k}^-+ \lambda_{u,k} / u_{k}^-},
\end{eqnarray}	
where $o( u^{H-1}\ln u)$ does not depend on $k,S$ and further
$$\limit{u} \sup_{\abs{k} \le N(u) } \eaE{\abs{\lambda -\lambda_{u,k}}}=0.$$
We analyse next  the uniform asymptotics of
$$p_k(u):=\pk{\int_{[0,\mathscr{T}_1+\epsilon]}\I(\sup_{s\in [0,S]}Z_{u,k}(s,t)>u_{k}^-)dt>x,
 \int_{[\mathscr{T}_2-\epsilon, \mathscr{T}_3+\epsilon]}\I(\sup_{s\in [0,S]}Z_{u,k}(s,t)>u_k^-+ \lambda_{u,k} / \dd{u_{k}^-})dt>y}$$
as $u\rw\IF$ with respect to $|k|\leq N(u)+1$.
\eeh{In order to apply Lemma \ref{the-weak-conv} in Appendix, we need to check conditions {\bf C1-C3} therein.  The first condition {\bf C1} follows immediately from \eqref{caj}. The second condition {\bf C2} is a consequence of \eqref{Cor}, while  ${\bf C3}$  follows from \eqref{caj}.  Consequently, using further \eqref{caj}, the application of the aforementioned lemma is justified and we obtain }
\BQN\label{uniform1}
\lim_{u\rw\IF}\sup_{|k|\leq N(u)+1}\left|\frac{p_k(u)}{\Psi(u_{k}^-)}-\mathcal{B}_H^{x,y}(\mathscr{T}_1+\epsilon; \eeh{\lambda},\mathscr{T}_2-\epsilon,\mathscr{T}_3+\epsilon)([0,S])\right|=0.
\EQN
Hence
\BQN\label{S1}
\Sigma_1^+(u)&\leq& \sum_{|k|\leq N(u)+1}\mathcal{B}_H^{x,y}(\mathscr{T}_1+\epsilon; \eeh{\lambda},\mathscr{T}_2-\epsilon,\mathscr{T}_3+\epsilon)([0,S])\Psi(u_{k}^-)\nonumber\\
&\leq& \mathcal{B}_H^{x,y}(\mathscr{T}_1+\epsilon;\eeh{\lambda},\mathscr{T}_2-\epsilon,\mathscr{T}_3+\epsilon)([0,S])\Psi(\uA)\sum_{|k|\leq N(u)+1}e^{-A^2u^{2(1-H)}\times \frac{(1- \epsilon)B}{2A}\Delta^2(u)\times (kS)^2}\nonumber\\
&\sim&\frac{\mathcal{B}_H^{x,y}(\mathscr{T}_1+\epsilon;\eeh{\lambda},\mathscr{T}_2-\epsilon,\mathscr{T}_3+\epsilon)([0,S])}{S}\frac{\sqrt{2}(AB)^{-1/2}(1-\epsilon)^{-1/2}}
{u^{1-H}\Delta(u)}\Psi(\uA)\int_\mathbb{R}e^{-t^2}dt\nonumber\\
&\sim&\frac{\mathcal{B}_H^{x,y}(\mathscr{T}_1;\eeh{\lambda},\mathscr{T}_2,\mathscr{T}_3)([0,S])}{S}\frac{\sqrt{2\pi}(AB)^{-1/2}}
{u^{1-H}\Delta(u)}\Psi(\uA),\quad u\rw\IF, \epsilon\rw 0.
\EQN
{\it \underline{Upper bound of $\Sigma\Sigma_1(u)$}}.
\eeh{Suppose for \eL{notational simplicity} that $\lambda=0$.}
\dd{Then $\uAA=\uA$ and}
\BQNY
\Sigma\Sigma_1(u)\leq \sum_{|k|\leq N(u)+1} (q_{k,1}(u)+q_{k,2}(u)),
\EQNY
where
\BQNY
q_{k,1}(u)&=&\pk{\sup_{t\in [0,T_3^*], s\in [kS,(k+1)S]}Z_u(s,t)>\dd{\uA}
 \sup_{t\in [0, T_3^*], s\in [(k+1)S, (k+1)S+\sqrt{S}]}Z_u(s,t)>\uA}\\
 &\leq& \pk{
 \sup_{t\in [0, T_3^*], s\in [(k+1)S, (k+1)S+\sqrt{S}]}\overline{Z}_u(s,t)>u_{k+1}^-},\\
 q_{k,2}(u)&=&\pk{\sup_{t\in [0,T_3^*], s\in [kS,(k+1)S]}Z_u(s,t)>\dd{\uA},
 \sup_{t\in [0, T_3^*], s\in [(k+1)S+\sqrt{S}, (k+2)S]}Z_u(s,t)>\uA}\\
 &\leq& \pk{\sup_{t\in [0,T_3^*], s\in [kS,(k+1)S]}\overline{Z}_u(s,t)>\dd{{u}_k^-},
 \sup_{t\in [0, T_3^*], s\in [(k+1)S+\sqrt{S}, (k+2)S ]}\overline{Z}_u(s,t)>u_{k+1}^-},
 \EQNY
 with
$$Z_u(s,t)=Z(t^*+\Delta(u)s,\Delta(u)t).$$
Analogously as in (\ref{uniform1}), we have that
\BQNY
\lim_{u\rw\IF}\sup_{|k|\leq N(u)+1}\left|\frac{\pk{
 \sup_{t\in [0, T_3^*], s\in [(k+1)S, (k+1)S+\sqrt{S}]}\overline{Z}_u(s,t)>u_{k+1}^-}}{\Psi(u_{k+1}^-)}-\overline{\mathcal{B}}_H^{0}(T_3^*)\overline{\mathcal{B}}_H^{0}(\sqrt{S})\right|=0.
\EQNY

Thus in view of  (\ref{S1})
\BQNY
\sum_{|k|\leq N(u)+1}q_{k,1}(u)&\leq& \sum_{|k|\leq N(u)+1}\overline{\mathcal{B}}_H^{0}(T_3^*)\overline{\mathcal{B}}_H^{0}(\sqrt{S})\Psi(u_{k+1}^-)\\
&\leq& \frac{\overline{\mathcal{B}}_H^{0}(T_3^*)\overline{\mathcal{B}}_H^{0}(\sqrt{S})}{S}\frac{\sqrt{2\pi}(AB)^{-1/2}}
{u^{1-H}\Delta(u)}\Psi(\uA),\quad u\rw\IF.
\EQNY
Additionally, in light of (\ref{Cor}) for $u$ sufficiently large
\BQN\label{Corbound}
{|s-s'|^{2H}+|t-t'|^{2H}}\leq 2(\uA)^2\left(1-Cor(\overline{Z}_u(s,t),\overline{Z}_u(s',t'))\right)\leq 4(|s-s'|^{2H}+|t-t'|^{2H})
\EQN
for all $|s|,|s'| \leq \frac{2\ln u}{u^{1-H}\Delta(u)}, t,t'\in [0,T_*]$.
Thus by  \cite[Cor 3.1]{KEP2016} there exist two positive constants $\mathcal{C}, \mathcal{C}_1$ such that for $u$ sufficiently large and $S>1$
\BQNY
q_{k,2}(u)\leq \mathcal{C}S^4e^{-\mathcal{C}_1S^{\frac{H}{2}}}\Psi(u_{k,k+1}^-), \quad
u_{k,l}^-=\min(u_{k}^-, u_{l}^-).
\EQNY
Hence
\BQNY
\sum_{|k|\leq N(u)+1}q_{k,2}(u)&\leq& \sum_{|k|\leq N(u)+1}\mathcal{C}S^4e^{-\mathcal{C}_1S^{\frac{H}{2}}}\Psi(u_{k,k+1}^-)\\
&\leq & \mathcal{C}S^3e^{-\mathcal{C}_1S^{\frac{H}{2}}}\frac{\sqrt{2\pi}(AB)^{-1/2}}
{u^{1-H}\Delta(u)}\Psi(\uA),\quad u\rw\IF.
\EQNY
Therefore we conclude that
\BQN\label{SS1}
\Sigma\Sigma_1(u)\leq \left(\frac{\overline{\mathcal{B}}_H^{0}(T_3^*)\overline{\mathcal{B}}_H^{0}(\sqrt{S})}{S}+\mathcal{C}S^3e^{-\mathcal{C}_1S^
{\frac{H}{2}}}\right)\frac{\sqrt{2\pi}(AB)^{-1/2}}
{u^{1-H}\Delta(u)}\Psi(\uA), \quad u\rw\IF.
\EQN
\eL{Note that if $\lambda\neq 0$, \kdd{the bound derived in \eqref{SS1} changes only by a multiplication by some constant}, 
which does not affect the \eHe{negligibility} of $\Sigma\Sigma_1(u)$}.

{\it \underline{Upper bound of $\Sigma\Sigma_2(u)$}}.
In light of (\ref{Corbound}) and applying \cite[Cor 3.1]{KEP2016}, we have that
\BQN\label{SS2}
\Sigma\Sigma_2(u)&\leq& \sum_{|k|,|l|\leq N(u)+1, l\geq k+2}\mathcal{C}S^4e^{-\mathcal{C}_1|l-k-1|^{H}S^{H}}\Psi(u_{k,l}^-)\nonumber\\
&\leq&\sum_{|k|\leq N(u)+1}\mathcal{C}S^4\Psi(u_{k}^-)\sum_{l=1}^\IF e^{-\mathcal{C}_1l^{H}S^{H}}\nonumber\\
&\leq &\sum_{|k|\leq N(u)+1}\mathcal{C}S^4e^{-\mathbb{Q}_6S^H}\Psi(u_{k}^-)\nonumber\\
&\leq &\mathcal{C}S^3e^{-\mathbb{Q}_6S^H}\frac{\sqrt{2\pi}(AB)^{-1/2}}
{u^{1-H}\Delta(u)}\Psi(\uA),\quad u\rw\IF.
\EQN
Consequently, as $u\rw\IF$
\BQN\label{Upper}
\pi_1(u)&\leq& \left(\frac{\mathcal{B}_H^{x,y}(\mathscr{T}_1;\eeh{\lambda},\mathscr{T}_2,\mathscr{T}_3)([0,S])}{S}+
\frac{\overline{\mathcal{B}}_H^{0}(T_3^*)\overline{\mathcal{B}}_H^{0}(\sqrt{S})}{S}
+\mathcal{C}S^3 [e^{-\mathcal{C}_1S^{\frac{H}{2}}}+e^{-\mathbb{Q}_6S^H}]\right)\nonumber\\
&&\times
\frac{\sqrt{2\pi}(AB)^{-1/2}}
{u^{1-H}\Delta(u)}\Psi(\uA).
\EQN

{\it $\diamond$ \underline{Lower bound of $\pi_1(u)$}}. \eeh{Again for notation simplicity we assume $\lambda=0$.}
Observe that
\BQNY
&&\frac{1}{\Delta(u)}\int_{[0,\tnu{1}]}\I(\sup_{|s-t-t^*|\leq (\ln u)/u^{1-H}}Z(s,t)>\uA)dt\\
&&\quad \geq \sum_{|k|\leq N(u)}\frac{1}{\Delta(u)}\int_{[0,\tnu{1}]}\I(\sup_{s\in I_k(u)}Z(s+t^*,t)>\uA)dt\\
&&\quad \quad -\sum_{|k|, |l|\leq N(u), k<l}\frac{1}{\Delta(u)}\int_{[0,T^*(u)/\ru]}\I(\sup_{s\in I_k(u)}Z(s+t^*,t)>\uA, \sup_{s\in I_l(u)}Z(s+t^*,t)>\uA)dt\\
&&\quad :=F_1(u)-F_2(u),\\
&&\frac{1}{\Delta(u)}\int_{[\tnu{2},\tnu{3}]}\I(\sup_{|s-t-t^*|\leq (\ln u)/u^{1-H}}Z(s,t)>\uA)dt\\
&&\quad \geq \sum_{|k|\leq N(u)}\frac{1}{\Delta(u)}\int_{[\tnu{2},\tnu{3}	]}\I(\sup_{s\in I_k(u)}Z(s+t^*,t)>\uA)dt\\
&&\quad \quad -\sum_{|k|, |l|\leq N(u), k<l}\frac{1}{\Delta(u)}\int_{[0,\tnuu]}\I(\sup_{s\in I_k(u)}Z(s+t^*,t)>\uA, \sup_{s\in I_l(u)}Z(s+t^*,t)>\uA)dt\\
&&\quad :=F_3(u)-F_2(u).
\EQNY
Hence, for $0<\epsilon<1$ \eHe{(write $s_\epsilon=(1+ \epsilon) s)$)}
\BQNY
\pi_1(u)&\geq& \pk{F_1(u)-F_2(u)>x, F_3(u)-F_2(u)>y}\\
&\geq& \pk{F_1(u)> x_\epsilon , F_3(u)>x_\epsilon, F_2(u)<\epsilon\min(x,y)}\\
&\geq& \pk{F_1(u)>x_\epsilon, F_3(u)>y_\epsilon}-\pk{F_2(u)\geq \epsilon\min(x,y)}.
\EQNY
Note that
\BQNY
 &&\pk{F_1(u)> x_\epsilon, F_3(u)> y_\epsilon}\\
 &&\geq \pk{\exists |k|\leq N(u):\ndu \int_{[0,\tnu{1}]}\I(\sup_{s\in I_k(u)}Z(s+t^*,t)>\uA)dt>x_\epsilon, F_3(u)>y_\epsilon  }\\
 &&\geq \sum_{|k|\leq N(u)}\pk{\ndu \int_{[0,\tnu{1}]}\I(\sup_{s\in I_k(u)}Z(s+t^*,t)>\uA)dt>x_\epsilon, F_3(u)>y _\epsilon }\\
 &&\ \ -\Sigma\Sigma_1(u)-\Sigma\Sigma_2(u)\\
 &&\geq \Sigma_1^-(u)-\Sigma\Sigma_1(u)-\Sigma\Sigma_2(u),
\EQNY
and
$$\pk{F_2(u)\geq \epsilon\min(x,y)}\leq \pk{F_2(u)>0}\leq \Sigma\Sigma_1(u)+\Sigma\Sigma_2(u),$$
where
\BQNY
\Sigma_1^-(u)=\sum_{k=-N(u)}^{N(u)}
 \pk{\int_{[0,\mathscr{T}_1]}\I(\sup_{s\in [0,S]}Z_{u,k}(s,t)>u_k^-)dt> x_\epsilon,
 \int_{[\mathscr{T}_2, \mathscr{T}_3]}\I(\sup_{s\in [0,S]}Z_{u,k}(s,t)>u_k^-)dt> y_\epsilon}.
 \EQNY
Hence
\BQNY
\pi_1(u)\geq \Sigma_1^-(u)-2\Sigma\Sigma_1(u)-2\Sigma\Sigma_2(u).
\EQNY
Analogously as in (\ref{S1}), it follows that
\BQNY
\Sigma_1^-(u)\sim \frac{\mathcal{B}_H^{x,y}(\mathscr{T}_1;\eeh{\lambda},\mathscr{T}_2,\mathscr{T}_3)([0,S])}{S}\frac{\sqrt{2\pi}(AB)^{-1/2}}
{u^{1-H}\Delta(u)}\Psi(\uA),\quad u\rw\IF, \epsilon\to 0,
\EQNY
which together with the upper bound of $\Sigma\Sigma_i, i=1,2$ leads to
\BQN\label{Lower}
\pi_1(u)
&\geq& \left(\frac{\mathcal{B}_H^{x,y}(\mathscr{T}_1;\eeh{\lambda},\mathscr{T}_2,\mathscr{T}_3)([0,S])}{S}-
\frac{2\overline{\mathcal{B}}_H^{0}(T_3^*)\overline{\mathcal{B}}_H^{0}(\sqrt{S})}{S}
-2\mathcal{C}S^3[ e^{-\mathcal{C}_1S^{\frac{H}{2}}}+e^{-\mathbb{Q}_6S^H}]\right)\nonumber\\
&&\times
\frac{\sqrt{2\pi}(AB)^{-1/2}}
{u^{1-H}\Delta(u)}\Psi(\uA),\quad u\rw\IF.
\EQN
Next  by \Cref{lemA0}, we have
$$\lim_{S\rw\IF}\frac{\mathcal{B}_H^{x,y}(\mathscr{T}_1;\eeh{\lambda},\mathscr{T}_2,\mathscr{T}_3)([0,S])}{S}=
\mathcal{B}_H^{x,y}(\mathscr{T}_1;\eeh{\lambda},\mathscr{T}_2,\mathscr{T}_3)\in (0,\IF), \quad \lim_{S\rw\IF}\frac{\overline{\mathcal{B}}_H^{0}(T_3^*)\overline{\mathcal{B}}_H^{0}(\sqrt{S})}{S}=0.$$
Thus letting $S\rw\IF$ in (\ref{Upper}) and (\ref{Lower}) yields
\BQNY
\pi_1(u)\sim \mathcal{B}_H^{x,y}(\mathscr{T}_1;\eeh{\lambda},\mathscr{T}_2,\mathscr{T}_3)\frac{\sqrt{2\pi}(AB)^{-1/2}}
{u^{1-H}\Delta(u)}\Psi(\uA),\quad u\rw\IF,
\EQNY
which combined with (\ref{pi2}) leads to
\BQNY
\lefteqn{
\pk{ \nvu  \int_{[T_2(u),T_3(u)]}\I(Q(t)>\ruu)dt>y, \nvu \int_{[0,T_1(u)]}\I(Q(t)>\ru)dt>x}}\\
&\sim&
\mathcal{B}_H^{x,y}(\mathscr{T}_1;\eeh{\lambda},\mathscr{T}_2,\mathscr{T}_3)\frac{\sqrt{2\pi}(AB)^{-1/2}}
{u^{1-H}\Delta(u)}\Psi(\uA),\quad u\rw\IF
\EQNY
establishing the proof. \QED
\COM{
	{\it $\diamond$ \underline{Existence of Berman-type constant}}.
Combing (\ref{Upper}) and (\ref{Lower}), we have that
\BQNY
&&\frac{\mathcal{B}_H^{x,y}(T_1,T_2,T_3)([0,S])}{S}+\frac{\mathcal{B}_H^{x,y}(T_1,T_2,T_3)([0,\sqrt{S}])}{S}
+\mathcal{C}S^3 [ e^{-\mathcal{C}_1S^{\frac{H}{2}}}+e^{-\mathbb{Q}_6S^H}]\\
&&\quad \geq \frac{\mathcal{B}_H^{x,y}(T_1,T_2,T_3)([0,S_1])}{S_1}-\frac{2\mathcal{B}_H^{x,y}(T_1,T_2,T_3)([0,\sqrt{S_1}])}{S_1}
-2\mathcal{C}S_1^3 [ e^{-\mathcal{C}_1S_1^{\frac{H}{2}}}+ e^{-\mathbb{Q}_6S_1^H}],
\EQNY
implying that
\BQNY
\limsup_{S\rw\IF}\frac{\mathcal{B}_H^{x,y}(T_1,T_2,T_3)([0,S])}{S}
=\liminf_{S\rw\IF}\frac{\mathcal{B}_H^{x,y}(T_1,T_2,T_3)([0,S])}{S}<\IF.
\EQNY
Next we show the positivity of the above limit. Replacing $ |s-t-t^*|\leq (\ln u)/u^{1-H}$ by $t^*+\bigcup_{|2k|\leq N(u)} I_{2k}(u)$ in the derivation of lower bound of $\pi_1(u)$, we have
\BQNY
\pi_1(u)&\geq& \sum_{|2k|\leq N(u)}\pk{\int_{[0,\tnu{1}]}\I(\sup_{s\in I_k(u)}Z(s+t^*,t)>\uA)dt>(1+\epsilon)x}-2\Sigma\Sigma_2(u)\\
&\geq& \left(\frac{\mathcal{B}_H^{x,y}(T_1,T_2,T_3)([0,S])}{2S}-2\mathcal{C}S^3e^{-\mathbb{Q}_6S^H}\right)
\frac{\sqrt{2\pi}(AB)^{-1/2}}
{u^{1-H}\Delta(u)}\Psi(u),\quad u\rw\IF,
\EQNY
which together with (\ref{Upper}) leads to
\BQNY
&&\frac{\mathcal{B}_H^{x,y}(T_1,T_2,T_3)([0,S])}{S}+\frac{\mathcal{B}_H^{x,y}(T_1,T_2,T_3)([0,\sqrt{S}])}{S}
+\mathcal{C}S^3e^{-\mathcal{C}_1S^{\frac{H}{2}}}+\mathcal{C}S^3e^{-\mathbb{Q}_6S^H}\\
&&\quad \geq \frac{\mathcal{B}_H^{x,y}(T_1,T_2,T_3)([0,S_1])}{2S_1}-2\mathcal{C}S_1^3e^{-\mathbb{Q}_6S_1^H}.
\EQNY
Hence for $S_1$ sufficiently large
$$\liminf_{S\rw\IF}\frac{\mathcal{B}_H^{x,y}(T_1,T_2,T_3)([0,S])}{S}\geq \frac{\mathcal{B}_H^{x,y}(T_1,T_2,T_3)([0,S_1])}{2S_1}-2\mathcal{C}S_1^3e^{-\mathbb{Q}_6S_1^H}>0.$$
Consequently, we conclude that
\BQNY
\mathcal{B}_H^{x,y}(T_1,T_2,T_3)=\lim_{S\rw\IF}\frac{\mathcal{B}_H^{x,y}(T_1,T_2,T_3)([0,S])}{S}\in (0,\IF).
\EQNY
Letting $S\rw\IF$ in (\ref{Upper}) and (\ref{Lower}) yields that
\BQNY
\pi_1(u)\sim \mathcal{B}_H^{x,y}(T_1,T_2,T_3)\frac{\sqrt{2\pi}(AB)^{-1/2}}
{u^{1-H}\Delta(u)}\Psi(u),\quad u\rw\IF,
\EQNY
which combined with (\ref{pi2}) establishes the claim.\\
}

\prooftheo{TH1} Clearly, for all  $x,y$ non-negative with $\ruuu= \ruu$
\BQN
\mathscr{P}^{x,y}_{{T}_1,{T}_2,{T}_3}(\ruuu,u)=\frac{\pk{\nvu\int_{[0,T_1(u)]}\I(Q(t)>u)dt>x, \nvu \int_{[T_2(u),T_3(u)]}\I(Q(t)>\ruuu)dt>y }}{\pk{\nvu\int_{[0,T_1(u)]}\I(Q(t)>u)dt>x}}.\nonumber
\EQN
The asymptotics of the denominator and the nominator are derived in Proposition \ref{asymptotics} and Proposition \ref{double.as}, respectively. Hence, using further \eqref{ciediamo} establishes the claim.
\QED

\subsection{Proof of Theorem \ref{TH3}}
{\underline{Case $a\in(-1,0)$.}}
Observe that
\BQNY
\lefteqn{
\pk{  \nvu \int_{[T_2(u),T_3(u)]}\I(Q(t)>(1+a)u)dt>y, \nvu \int_{[0,T_1(u)]}\I(Q(t)>u)dt>x}}\\
&=&
\pk{  \nvu \int_{[0,T_1(u)]}\I(Q(t)>u)dt>x}\\
&&-
\pk{ \nvu  \int_{[T_2(u),T_3(u)]}\I(Q(t)\le(1+a)u)dt>T_3(u)-T_2(u)-y, \nvu  \int_{[0,T_1(u)]}\I(Q(t)>u)dt>x}\\
&=:& P_1(u)-P_2(u).
\EQNY
Next, recalling that $T^*(u)=\max(T_1(u), T_2(u), T_3(u))$ and using that $T^*(u)\sim C u^{(2H-1)/H}$ as $u\to\infty$ for some $C>0$, we obtain
\BQNY
P_2(u)&\le&\pk{\inf_{t\in [T_2(u),T_3(u)]} Q(t)\leq (1+a)u, \sup_{t\in [0,T_1(u)]}Q(t)>u} \\
&\le&\pk{\text{there exist}~t, s \in [0,T^*(u)], Q(t)-Q(s)\geq -au }\\
&\le&
\rd{\pk{\sup_{0\le t\le s \le T^*(u) } (B_H(t)-B_H(s)-c(t-s))> -au}}\\
&\le&
\pk{\sup_{0\le t\le s \le 1 } {T^*}^H(u)(B_H(t)-B_H(s))> -au}\\
&\le&
{\rm C_1}e^{-{\rm C_2}u^{4-4H}}
\EQNY
for some $C_1,C_2>0$,
where the third inequality is because of \eqref{def.Q} and  the last inequality above is due to Borell-TIS inequality.
Hence, in view of Proposition \ref{asymptotics},
$P_2(u)=o(P_1(u))$ as $u\to\infty$, which
  leads to
\[
\pk{ \nvu  \int_{[T_2(u),T_3(u)]}\I(Q(t)>(1+a)u)dt>y \Big{|} \nvu  \int_{[0,T_1(u)]}\I(Q(t)>u)dt>x}\sim 1
\]
as $u\to\infty$.\\
{\underline{Case $a>0$.}}
First, we consider the asymptotic upper bound. We note that
\BQNY
\mathscr{P}^{x,y}_{{T}_1,{T}_2,{T}_3}((1+a)u,u)
\le
\frac{\pk{ \nvu \int_{[T_2(u),T_3(u)]}\I(Q(t)>(1+a)u)dt>y}}
     {\pk{\nvu\int_{[0,T_1(u)]}\I(Q(t)>u)dt>x}}
\EQNY
and
\BQNY
\lefteqn{
\pk{ \nvu \int_{[T_2(u),T_3(u)]}\I(Q(t)>(1+a)u)dt>y}}\\
&=&
\pk{ \frac 1 {v((1+a)u)} \int_{[T_2(u),T_3(u)]}\I(Q(t)>(1+a)u)dt> (1+a)^{(1-2H)/H}y},
\EQNY
where, by {\bf (T)} we have
$$\lim_{u\to\infty}\frac{T_i(u)}{v((1+a)u)}=\mathscr{T}_i(1+a)^{(1-2H)/H}, \quad i=1,2.
$$
Consequently,  by the stationarity of $Q(t), t\ge0$ and Proposition \ref{asymptotics}, with
$\tilde{a}=(1+a)^{(1-2H)/H}$ we obtain
\BQNY
\lefteqn{
\pk{ \nvu \int_{[T_2(u),T_3(u)]}\I(Q(t)>(1+a)u)dt>y}}\\
&=&
\pk{ \nvu \int_{[0, T_3(u)-T_2(u)]}\I(Q(t)>(1+a)u)dt>y}\\
&\sim&
\mathcal{H}_{2H}\overline{\mathcal{B}}_H^{{\tilde{a}}y}((\mathscr{T}_3-\mathscr{T}_2){\tilde{a}})\frac{\sqrt{2\pi}(AB)^{-1/2}}
{(1+a)^{1-H}u^{1-H}\Delta((1+a)u)}\Psi(A((1+a)u)^{1-H})
\EQNY
as $u\to\infty$.
Hence
\BQNY
\limsup_{u\to\infty}
\frac{\mathscr{P}^{x,y}_{{T}_1,{T}_2,{T}_3}((1+a)u,u)}
{\exp\left(-\frac{A^2\left((1+a)^{2-2H}-1\right)}{2}u^{2-2H}\right)}
\le
\tilde{a}^{1-H}\frac{\overline{\mathcal{B}}_H^{\tilde{a}y}(\tilde{a}(\mathscr{T}_3-\mathscr{T}_2))}
{\overline{\mathcal{B}}_H^{x}(\mathscr{T}_1)}.
\EQNY

For the proof of the asymptotic lower bound we have
\BQNY
\lefteqn{
\mathscr{P}^{x,y}_{{T}_1,{T}_2,{T}_3}((1+a)u,u)}\\
&\ge&
\frac{\pk{ \nvu \int_{[T_2(u),T_3(u)]}\I(Q(t)>(1+a)u)dt>y,  \nvu\int_{[0,T_1(u)]}\I(Q(t)>(1+a)u)dt>x}}
     {\pk{ \nvu\int_{[0,T_1(u)]}\I(Q(t)>u)dt>x}}.
\EQNY
Then, following the same line of arguments as for the asymptotic upper bound,
by Proposition \ref{double.as}
we obtain
\BQNY
\liminf_{u\to\infty}
\frac{\mathscr{P}^{x,y}_{{T}_1,{T}_2,{T}_3}((1+a)u,u)}
{\exp\left(-\frac{A^2\left((1+a)^{2-2H}-1\right)}{2}u^{2-2H}\right)}
\ge
\tilde{a}^{1-H}
\frac{
\overline{\mathcal{B}}_H^{\tilde{a}x,\tilde{a}y}(\tilde{a}\mathscr{T}_1;\eeh{0},\tilde{a} \mathscr{T}_2,\tilde{a}\mathscr{T}_3)}
{\overline{\mathcal{B}}_H^{x}(\mathscr{T}_1)}
.
\EQNY
\QED

\section{Appendix}

In this Section we present a lemma that plays a crucial lemma for proof of \Cref{double.as}. Consider next
$$\xi_{u,j}(s,t), \quad (s, t)\in \KK=[0,S]\times [0,T], \quad j\in S_u$$ a family of centered Gaussian random fields with continuous sample paths and unit variance, \eL{where $S_u$ is \eHe{a countable}  index set}.  For $S>0$, $0<b_1, b_2, b_3\leq \eeh{T}$, $b_1>x\geq 0$ and $b_3-b_2>y\geq 0$, we are interested in the uniform asymptotics of
$$ p_{u,j}(S;\eeh{\lambda_{u,j}})=\pk{ \int_{[0,b_1]} \mathbb{I}\LT(\sup_{s\in[0, S]}\xi_{u,j}(s,t)>g_{u,j}\RT) dt >x, \int_{[b_2,b_3]} \mathbb{I}\LT(\sup_{s\in[0, S]}\xi_{u,j}(s,t)>g_{u,j}+ \lambda_{u,j}/g_{u,j}\RT) dt >y}$$
 with respect to $j\in S_u$, as $u\to\infty$, where \eHe{$g_{u,j}$'s  and $\lambda_{u,j}$'s}  are given constants depending on $u$ and $j$.

 Suppose next that $S_u$'s are finite index. The following assumption\eHe{s} will be imposed in the lemma below:
 \begin{itemize}
 	\item[{\bf C1}:]   $g_{u,j}, j\in S_u, u>$  \eHe{are constants }  satisfying
 	\BQNY
 	\lim_{u\to\IF}\inf_{j\in S_u}g_{u,j}=\IF.
 	\EQNY
 	\item[{\bf C2}:] There exists $\alpha\in (0,2]$ such that
 \BQNY
\lim_{u\rw\IF}\sup_{j\in S_u}\sup_{(s,t)\neq (s',t'), (s,t), (s',t')\in E} \left|g_{u,j}^2\frac{1-Corr(\xi_{u,j}(s,t), \xi_{u,j}(s',t'))}{|s-s'|^{\alpha}+|t-t'|^{\alpha}}-1\right|=0.
 \EQNY
  	\item[{\bf C3}:]
  	\rd{The sequence $\lambda_{u,j}$ is such that
 	\BQNY
\lim_{u\to\IF}\sup_{j\in S_u}|\lambda_{u,j}-\lambda|=0
\EQNY
for some $\lambda\in\R$.
}
 \end{itemize}

\eHe{We state next} a modification of \cite[Lem 4.1]{DHLM23}.
\BEL\label{the-weak-conv}
Let $\{\xi_{u,j}(s,t), (s,t)\in E ,j\in S_u\}$ be a family of centered Gaussian random fields defined as above. If {\bf C1-C3} holds, then for all $S>0$, $0<b_1, b_2, b_3\leq b$, $b_1>x\geq 0$ and $b_3-b_2>y\geq 0$ we have
\BQN\label{con-uni-con}
\lim_{u\to\IF}\sup_{j\in S_u} \ABs{ \frac{p_{u,j}(S; \eeh{\lambda_{u,j}} )} {\Psi(g_{u,j})} -\mathcal{B}_{\alpha/2}^{x,y} (b_1;\eeh{\lambda},b_2,b_3)([0, S]) } = 0.
\EQN
\EEL

\prooflem{the-weak-conv}
The proof of Lemma \ref{the-weak-conv} follows by similar argumentation as
given in the proof of \cite[Lemma 4.1]{DHLM23}. For completeness, we present details of the
main steps of the argumentation.

Let
$$ \chi_{u,j}(s,t) := g_{u,j}(\xi_{u,j}(s,t)-\rho_{u,j}(s,t)\xi_{u,j}(0,0)),\quad (s,t)\in E,
$$
and
$$f_{u,j}(s,t,w):=w\rho_{u,j}(s,t)-g^2_{u,j}\LT( 1-\rho_{u,j}(s,t) \RT), \quad (s,t)\in E, w\in\R,$$
where $\rho_{u,j}(s,t)=Cov(\xi_{u,j}(s,t), \xi_{u,j}(0,0))$.
Conditioning on $\xi_{u,j}(0,0)$ and using the fact that $\xi_{u,j}(0,0)$ and $\xi_{u,j}(s,t)-\ruk(s,t)\xi_{u,j}(0,0)$ are mutually independent, we obtain
\BQNY
&& p_{u,j}(S;\lambda_{u,j})\\
&&=\frac{e^{-g^2_{u,j}/2}} { \sqrt{2\pi}g_{u,j}}
\int_{\R} \exp\LT(-w-\frac{w^2}{2g^2_{u,j}}\RT) \pk{ \int_{0}^{b_1} \I\left(\sup_{s\in [0,S]}\left(g_{u,j}(\xi_{u,j}(s,t)-g_{u,j}) \right)>0\right) dt > x,\right.\\
&&\quad\left. \int_{b_2}^{b_3} \I\left(\sup_{s\in [0,S]}\left(g_{u,j}(\xi_{u,j}(s,t)-g_{u,j})-\eeh{\lambda}_{u,j}\right)>0\right) dt > y \Big{|} \xi_{u,j}(0,0)=g_{u,j}+wg_{u,j}^{-1}} dw\\
&&= \frac{e^{-g^2_{u,j}/2}} { \sqrt{2\pi}g_{u,j}}
\int_{\R} \exp\LT(-w-\frac{w^2}{2g^2_{u,j}}\RT) \pk{\int_{0}^{b_1} \I\left(\sup_{s\in [0,S]}\left(\chi_{u,j}(s,t)+f_{u,j}(s,t,w)\right)>0\right) dt > x,\right.\\
&&\quad\left. \int_{b_2}^{b_3} \I\left( \sup_{s\in [0,S]}\left(\chi_{u,j}(s,t)+f_{u,j}(s,t,w)- \eeh{\lambda}_{u,j}\right)> 0\right ) dt > y } dw\\
&&:=\frac{e^{-g^2_{u,j}/2}} { \sqrt{2\pi}g_{u,j}} \int_{\R} \exp\LT(-w-\frac{w^2}{2g^2_{u,j}}\RT) \mathcal{I}_{u,j}(w; x,y) dw,
\EQNY
where
\BQNY\mathcal{I}_{u,j}(w;x,y)&=&\pk{\int_{0}^{b_1} \I\left(\sup_{s\in [0,S]}\left(\chi_{u,j}(s,t)+f_{u,j}(s,t,w)\right)>0\right) dt > x,\right.\\
&&\quad \left. \int_{b_2}^{b_3} \I\left(\sup_{s\in [0,S]}\left(\chi_{u,j}(s,t)+f_{u,j}(s,t,w)- \eeh{\lambda}_{u,j}\right)>0\right) dt > y }.
\EQNY
Noting that
$$\lim_{u\rw\IF}\sup_{j\in S_u}\left|\frac{\frac{e^{-g^2_{u,j}/2}} { \sqrt{2\pi}g_{u,j}}}{\Psi(g_{u,j})}-1\right|=0
$$
and for any $M>0$
$$\lim_{u\to\infty}\inf_{|w|\leq M}e^{-\frac{w^2}{2g^2_{u,j}}}=1
$$
we can establish the claim if we show that
\BQN\label{convergence}
\lim_{u\to\IF}\sup_{j\in S_u} \ABs{ \int_{\R} \exp\LT(-w\RT) \mathcal{I}_{u,j}(w,x,y) dw - \mathcal{B}_{\alpha/2}^{x,y} (b_1;\eeh{\lambda},b_2,b_3)([0, S])} = 0.
\EQN

{\bf \underline{Weak convergence}}.
We next show the weak convergence of  $\{\chi_{u,j}(s,t) + f_{u,j}(s,t,w), (s,t)\in E\}$ as $u\to\infty$.
By {\bf C1} and {\bf C2} we have,  for $(s,t), (s',t')\in E$, as $u\to\infty$, uniformly with respect to $j\in S_u$
\BQNY Var(\chi_{u,j}(s,t)-\chi_{u,j}(s',t')) &=& g^2_{u,j}\LT( \E{\xi_{u,j}(s,t)-\xi_{u,j}(s',t')}^2 - \LT( \rho_{\xi_{u,j}}(s,t)-\rho_{\xi_{u,j}}(s',t') \RT)^2 \RT) \\
&&\quad \rightarrow 2Var(\zeta(s,t)-\zeta(s',t')),
\EQNY
where $\zeta(s,t)=B_{\alpha/2}(s)+B'_{\alpha/2}(t), ~(s,t)\in E$ with $B$ and $B'$ being independent fBm's.
This implies that the finite-dimensional distributions of $\{\chi_{u,j}(s,t), (s,t)\in E\}$ weakly converge to that of $\{\sqrt{2}\zeta(s,t), (s,t)\in E\}$ as $u\to\infty$  uniformly with respect to $j\in S_u$.
Moreover, it follows from {\bf C2} that, for $u$ sufficiently large
$$Var(\chi_{u,j}(s,t)-\chi_{u,j}(s',t'))\leq g^2_{u,j} \E{\xi_{u,j}(s,t)-\xi_{u,j}(s',t')}^2\leq 4(|s-s_1|^{\alpha}+|t-t_1|^{\alpha}),~(s,t),(s_1,t_1)\in E.$$
This implies that uniform tightness of $\{\chi_{u,j}(s,t), (s,t)\in E\}$ for large $u$ with respect to $j\in S_u$. Hence
$\{\chi_{u,j}(s,t), (s,t)\in E\}$ weakly converges to $\{\sqrt{2}\zeta(s,t), (s,t)\in E\}$ as $u\to\infty$ uniformly with respect to $j\in S_u$.
Additionally, by
{\bf C1}-{\bf C2},
$\{f_{u,j}(s,t,w), (s,t)\in E\}$ converges to $\{w-|s|^{\alpha}-|t|^{\alpha}, (s,t)\in E\}$ uniformly with respect to $j\in S_u$.
Therefore, we conclude that as $u\to\infty$,
$\{\chi_{u,j}(s,t) + f_{u,j}(s,t,w), (s,t)\in E\}$ weakly converges to $\{\sqrt{2}\zeta(s,t)+w-|s|^{\alpha}-|t|^{\alpha}, (s,t)\in E\}$ uniformly with respect to $j\in S_u$.
Then continuous mapping theorem implies that
$$\{z_{u,j}(t,w)=\sup_{s\in [0,S]}\left(\chi_{u,j}(s,t) + f_{u,j}(s,t,w)\right), t\in [0, b]\}$$ weakly converges to $$\{z(t)+w=\sup_{s\in [0,S]}\left(\sqrt{2}\zeta(s,t)+w-|s|^{\alpha}-|t|^{\alpha}\right), t\in [0,b]\}$$ uniformly with respect to $j\in S_u$ for each $w\in\mathbb{R}$.\\
\eeh{Repeating the arguments, in view of {\bf C3} the same convergence holds for
	$\chi_{u,j}(s,t) + f_{u,j}(s,t,w)+ \lambda_{u,j}$.}\\
In order to show the weak convergence of
$$\left(\int_{0}^{b_1}\mathbb{I}(z_{u,j}(t,w)>0)dt,  \int_{b_2}^{b_3}\mathbb{I}(z_{u,j}(t,w)- \lambda_{u,j}>0)dt\right),\quad u\to \IF$$
 we have to prove that
$$\left(\int_{0}^{b_1}\mathbb{I}(f(t)>0)dt, \int_{b_2}^{b_3}\mathbb{I}(f(t)>\eeh{\lambda})dt\right)$$  is a continuous functional from $C([0,b_1]\cup[b_2,b_3])$ to $\mathbb{R}^2$  except a zero probability subset of  $C([0,b_1]\cup[b_2,b_3])$ under the probability induced by $\{z(t)+w, t\in [0,b_1]\cup [b_2,b_3]\}$. The idea of the proof follows from Lemma 4.2 of \cite{berman1973excursions}. Observe that the discontinuity set is
$$E^*=\left\{f\in C([0,b_1]\cup[b_2,b_3]): \int_{[0,b_1]}\mathbb{I}(f(t)=0)    dt>0~\text{or}~\int_{[b_1,b_2]}\mathbb{I}(f(t)=\lambda)    dt>0\right\}.$$
Note that for any $c\in\R$
$$\int_\mathbb{R} \mathbb{E}\left(\int_{[0,b_1]\cup[b_2,b_3]}\mathbb{I}(z(t)+w=c)dt\right)dw=\int_{[0,b_1]\cup[b_2,b_3]}\int_\mathbb{R}\pk{z(t)+w=c}dwdt=0.$$
Hence
$E^*$  has probability zero under the probability induced by $\{z(t)+w, t\in [0,b_1]\cup[b_2,b_3]\}$ for a.e. $w\in\mathbb{R}$.
Application of the continuous mapping theorem yields that $$\left(\int_{0}^{b_1}\mathbb{I}(z_{u,j}(t,w)>0)dt, \int_{b_2}^{b_3}\mathbb{I}(z_{u,j}(t,w)>\lambda)dt\right)$$ weakly converges to $$\left(\int_{0}^{b_1}\mathbb{I}(z(t)+w> 0 )dt, \int_{b_2}^{b_3}\mathbb{I}(z(t)+w>\lambda)dt\right)$$ as $u\to\infty$,  uniformly with respect to $j\in S_u$ for a.e. $w\in\mathbb{R}$.

{\bf \underline{Convergence on continuous points}}.
Let
$$\mathcal{I}(w;x,y):=\pk{\int_{0}^{b_1}\mathbb{I}(z(t)+w>0)dt>x, \int_{b_2}^{b_3}\mathbb{I}(z(t)+w>\lambda)dt>y}.$$
Using  similar arguments as  in the proof of
\cite[Thm 1.3.1]{Berman92}, we show that (\ref{convergence}) holds for continuity points $(x,y)$ with $x, y>0$, i.e.,
$$ \lim_{\vp\to0}\int_{\R} \LT(\mathcal{I}(w;x+\vp, y+\vp) -\mathcal{I}(w;x-\vp, y-\vp)\RT)e^{-w} dw=0.$$
Note that for all  $x, y>0$
\BQN\label{PiterbargTh}\mathcal{I}(w;x,y)&\leq& \pk{\sup_{(s,t)\in E} \sqrt{2}\zeta(s,t)-|s|^{\alpha}-|t|^{\alpha}>-w}\nonumber\\
&\leq& \pk{\sup_{(s,t)\in E} \sqrt{2}\zeta(s,t)>-w+C}\nonumber\\
&
\leq& C_1e^{-Cw^2},~ w<-M
\EQN
for $M$ sufficiently large, where $C, C_1$ \eeh{are positive constants} and in the last inequality, we \eeh{used the} Piterbarg inequality \cite[Thm 8.1]{Pit96}. Hence the dominated convergence theorem gives
$$\int_{\R} \LT(\mathcal{I}(w;x+, y+) -\mathcal{I}(w;x-,y-)\RT)e^{-w} dw=0.$$
This implies that if $(x,y)$ is a continuity point, then $\mathcal{I}(w;)$ is continuous at $(x,y)$ for a.e. $w\in \mathbb{R}$. Hence if $(x,y)$ is a continuity point, then
\BQN\label{con1}\lim_{u\to\infty}\sup_{j\in S_u}\left|\mathcal{I}_{u,j}(w;x,y)-\mathcal{I}(w;x,y)\right|=0,~\text{for a.e.} ~w\in\mathbb{R}.
\EQN
\eeh{Applying again the} Piterbarg inequality, analogously as in (\ref{PiterbargTh}), we obtain
\BQN\label{Pit1}\sup_{j\in S_u}\mathcal{I}_{u,j}(w;x,y)\leq  C_1e^{-Cw^2},~ w<-M
\EQN
for $M$ and $u$ sufficiently large.
\eeh{Consequently,} in view of  (\ref{PiterbargTh}), (\ref{con1}) and (\ref{Pit1}), \eeh{the} dominated convergence theorem establishes (\ref{convergence}).

{\bf \underline{Continuity of $\mathcal{B}_{\alpha/2}^{x,y} (b_1; \eeh{\lambda},b_2,b_3)([0, S])$}}.
Clearly, $\mathcal{B}_{\alpha/2}^{x,y} (b_1;\lambda, b_2,b_3)([0, S])$ is right-continuous at $(x,y)=(0,0)$. We next focus on its continuity over $\left([0, b_1)\times [0,b_3-b_2)\right)\setminus \{(0,0)\}$. To show $\mathcal{B}_{\alpha/2}^{x,y} (b_1;\lambda, b_2,b_3)([0, S])$ is continuous at $(x,y)\in (0,b_1)\times (0,b_3-b_2)$, it suffices to prove that
$$\int_\mathbb{R}e^{-w}
\left(\pk{ \int_{0}^{b_1}\mathbb{I}(z(t)+w>0)dt=x}+\pk{\int_{b_2}^{b_3}\mathbb{I}(z(t)+w>\lambda)dt=y}\right)dw=0.$$
Denote $A_w=\{z_\kappa(t): \int_{0}^{b_1}\mathbb{I}(z_\kappa(t)+w>0)dt=x\}$, where $z_\kappa(t)=z(t)(\kappa)$ with $\kappa\in\Omega$ the sample space. In light of the continuity of $z_\kappa(t)$, if $\int_{0}^{b_1}\mathbb{I}(z_\kappa(t)+w>0)dt=x$ for $x\in (0,b_1)$ and $w'>w$, then
$$\int_{0}^{b_1}\mathbb{I}(z_\kappa(t)+w'>0)dt>x.$$
Hence $A_w\cap A_{w'}=\emptyset$ if $w\neq w'$.
Noting that the continuity of $z(s)$ guarantees the measurability of $A_w$, and
$$\sup_{\Lambda\subset \mathbb{R}, \#\Lambda<\infty}\sum_{w\in \Lambda}\pk{A_w}\leq 1,$$
where  $\#\Lambda$ stands for  the cardinality of the set $\Lambda$.\\
\eeh{Note in passing the important fact that $\mathbb{P}$-measurability of $A_w$ is a consequence of the Fubini-Tonelli theorem.}
Next, it follows that
$$\{w: w\in \mathbb{R} \quad \text{such that} \quad \pk{A_w}>0\}$$
is a countable set, which implies that for $x\in (0,b_1)$
$$\int_\R  \pk{A_w} e^{-w} dw=0. $$
Using similar argument, we can show
$$\int_\mathbb{R}e^{-w}
\pk{\int_{b_2}^{b_3}\mathbb{I}(z(t)+w>\lambda)dt=y}dw=0.$$
Therefore, we conclude that
$\mathcal{B}_{\alpha/2}^{x,y} (b_1;\lambda,b_2,b_3)([0, S])$ is continuous at $(x,y)\in (0,b_1)\times (0,b_3-b_2)$. Analogously, we can show the continuity on $\{0\}\times (0,b_3-b_2)$ and $(0,b_1)\times\{0\}$.
This completes the proof. \QED

{\bf Acknowledgement}: K. D\c ebicki
was partially supported by NCN Grant No 2018/31/B/ST1/00370 (2019-2024).
Financial support from the Swiss National Science Foundation Grant 200021-196888 is also kindly acknowledged.

 \bibliographystyle{ieeetr}

 \bibliography{queue3c}

\end{document}